\documentclass[reqno,12pt]{amsart}
\usepackage{amsfonts,amsmath,amsthm, url}
\usepackage{amssymb}
\usepackage{graphics}
\usepackage{graphicx}
\usepackage{epsfig,color,subfigure}
   \newtheorem{thm}{Theorem}
   \newtheorem{prop}{Proposition}
   \newtheorem{lem}{Lemma}[section]

   \newtheorem{rem}{Remark}[section]

\newcommand{\N}{\mathbb{N}}

\newcommand{\R}{\mathbb{R}}
\newcommand{\C}{\mathbb{C}}
\newcommand{\E}{\mathbb{E}}
\newcommand{\prob}{\mathbb{P}}
\newcommand{\p}{\partial}

\newcommand{\eps}{\varepsilon}

\newcommand{\ds}{\displaystyle}

\def\un{{\mathrm{1~\hspace{-1.4ex}l}}}

\makeatletter
  
  \@addtoreset{equation}{section}
\makeatother
\begin{document}
\title[Long time behavior of GP equation at positive temperature]
{Long time behavior of Gross-Pitaevskii equation at positive temperature}
\author[Anne DE BOUARD]{Anne DE BOUARD$^{\scriptsize 1}$}
\author[Arnaud Debussche]{Arnaud Debussche$^{\scriptsize 2}$}
\author[Reika FUKUIZUMI]{Reika FUKUIZUMI$^{\scriptsize 3}$}



\maketitle

\begin{center} \small
$^1$ CMAP,  Ecole Polytechnique, CNRS, Universit\'e Paris-Saclay, \\
91128 Palaiseau, France; \\
\email{debouard@cmap.polytechnique.fr}
\end{center}

\begin{center} \small
$^2$ IRMAR, ENS Rennes, CNRS, UBL. \\ 
av Robert Schuman, F-35170 Bruz, France; \\
\email{arnaud.debussche@ens-rennes.fr}
\end{center}

\begin{center} \small
$^3$ Research Center for Pure and Applied Mathematics, \\
Graduate School of Information Sciences, Tohoku University,\\
Sendai 980-8579, Japan; \\
\email{fukuizumi@math.is.tohoku.ac.jp}
\end{center}

\vskip 0.1 in
\noindent
{\small 
{\bf Abstract}. The stochastic Gross-Pitaevskii equation is used as a model 
to describe Bose-Einstein condensation at positive temperature. The equation 
is a complex Ginzburg Landau equation with a trapping potential and 
an additive space-time white noise. Two important questions for this system 
are the global existence of solutions in the support of the Gibbs measure, and the convergence 
of those solutions to the equilibrium for large time. In this paper, we give a proof of 
these two results in one space dimension. In order to prove the convergence to equilibrium, 
we use the associated purely dissipative equation as an auxiliary equation, for which the convergence may
be obtained using standard techniques. Global existence is obtained for all initial data, and not almost surely with respect
to the invariant measure.
\vskip 0.1 in

\section{Introduction} 

In this paper, we will present a mathematical analysis of a model related  to discussions on the Gibbs equilibrium   
in the papers \cite{bbdbg,kc}. In those papers, the authors  consider the dynamics 
of the wave function in Bose-Einstein condensation 
near the critical temperature $T_c$, using the so called (projected) stochastic Gross-Pitaevskii equation :
\begin{equation} \label{eq:SGPE}
d\psi=\mathcal{P}\Big\{-\frac{i}{\hbar}L_{GP}\psi dt +\frac{G(x)}{k_B T}(\nu-L_{GP})\psi dt +dW_G(x,t)\Big\}
\end{equation}
where 
\begin{eqnarray*}
L_{GP}=-\frac{\hbar^2}{2m} \nabla^2 +V(x)+g|\psi|^2,
\quad \langle dW_G^*(s,y), dW_G(t,x)\rangle= 2G(x)\delta_{t-s} \delta_{x-y} \, dt,
\end{eqnarray*}
whose derivation from the quantum mechanics system may be found in, e.g., \cite{ds,gd}. Here, 
$m$ is the mass of an atom, $V(x)$ is the trapping potential, $\nu$ is the chemical potential,  
$g$ characterizes the strength of atomic interactions related to the s-wave scattering
length. The second and third terms in the right hand side of (\ref{eq:SGPE})
represent growth processes, i.e., collisions that transfer atoms from the
thermal cloud to the classical field and vice versa. 
The form of $G(x)$ may be determined from kinetic theory,  
and is often taken as a constant, and $dW_G$ is the complex-valued Gaussian noise associated 
with the condensate growth. Lastly, $\mathcal{P}$ is a projection which 
restricts the dynamics to the low-energy region defined 
by the harmonic oscillator modes, or Fourier modes, depending on the situation. 
At zero temperature $T=0$, the statistics of the atoms is well represented 
by a single condensate wave function, and the standard Gross-Pitaevskii 
equation describes the coherent evolution of the wave function in a quite good manner 
since all spontaneous and incoherent processes, like the effect of thermal cloud, 
may be neglected. The effects of such incoherent elements 
are implemented by adding a dissipation and a noise to the standard Gross Pitaevskii equation as above.
This model was numerically shown in \cite{w}, 
to be well suited for modeling the dynamics of evaporative cooling and vortex formation during BEC.
The authors in \cite{bbdbg} propose a more detailed high-energy cut-off stochastic Gross-Pitaevskii 
equation (which they call SPGPE), 
where in their model a rotation term and a multiplicative noise are also considered.

The authors in \cite{kc} use (\ref{eq:SGPE}) with $V \equiv 0$ and periodic boundary condition in 2d or 3d 
for the purpose of classification of the phase transition in the system 
by means of well-known universality class argument in thermal equilibrium. 
The phase transition accompanied with a symmetry breaking leads 
to the formation of a topological defect, e.g., strings, vortices, 
and the paper \cite{kc} tries to find some common properties 
between the strings at the beginning of the Universe and the vortices in Bose-Einstein condensates 
by studying the thermal (Gibbs) equilibrium and the stochastic dynamic evolution. 
\vspace{3mm}

Our purpose in the present paper is the study of the stochastic partial differential equation (\ref{eq:SGPE}) 
with $\mathcal{P}=Id$ and $G(x)$ constant (which we assume from now on). Using the mathematical construction of the Gibbs equilibrium measure 
in dimension one in space studied in \cite{btt}, together with the strong Feller property of the stochastic evolution
equation, we will prove the global existence of solutions for the infinite dimensional system, and
show that those solutions  converge exponentially to the Gibbs equilibrium. Note that this implies, for the finite
dimensional system \eqref{eq:SGPE}, the exponential convergence with a rate which is uniform in term of 
the number of modes taken into account in the projection $\mathcal{P}$.
Several related mathematical studies can be found; for example, in \cite{bs} the author dealt 
with the equation (\ref{eq:SGPE}) on a bounded domain $D$ with $V\equiv 0$ (and without the projection $\mathcal{P}$), 
in both cases of defocusing and focusing nonlinearities. 
It was shown there the existence and uniqueness of the invariant measure on $L^p(D)$ for any $p \ge 2$, 
but the convergence to the Gibbs measure was not discussed. Recently, the authors in \cite{cfl} 
considered an equation close to (\ref{eq:SGPE}) with a focusing nonlinearity (i.e. $g<0$), $V\equiv 0$ 
and a regular noise, but with a modified linear operator, on the $1d$ torus. They proved the existence and the invariance 
of the grand-canonical Gibbs measure, and the exponential decay to equilibrium. 
The additional term added in the dynamical equation in \cite{cfl} may be viewed 
as a restoring term. The techniques we develop in the present paper allow to simplify the proof of the convergence 
to equilibrium in this case, and to treat the case of the space-time white noise (see Remark \ref{rem7.2}).

\vspace{3mm}

We will consider the stochastic Gross-Pitaevskii equation (\ref{eq:SGPE}) on $\R$, 
a defocusing nonlinearity ($g>0$), and a harmonic potential $V(x)=|x|^2$.  
The invariance of the Gibbs measure, global existence of solutions for all initial data,
and the exponential decay to equilibrium 
will be proved. The support of the Gibbs measure 
in our case will be on a Banach space, 
and this fact requires a bit more complicated justification than the papers \cite{bs, cfl} where 
the basic space was the Hilbert space $L^2$. Also, the paper \cite{bs} assumed that the dissipation  
was not too small to obtain a globally defined strong solution, whereas  
our result covers any coefficient size of the dissipation.
The noise we will consider is a space-time white noise, 
and our proof for the convergence to the Gibbs equilibrium, together with the global existence
for all initial data thanks to the strong Feller property
are the first results for this kind of stochastic Ginzburg-Landau equations, 
as far as we know.  Also, due to the presence of the space-time white noise, the (more physical) case of space
dimension two or three requires the use of renormalization, and of much more involved arguments, and will be 
the object of future work. The focusing case also needs some investigations, as the techniques of 
\cite{cfl} would lead to a trivial measure in the presence of the quadratic potential (see again Remark \ref{rem7.2}).
\vspace{3mm}

Lastly, let us introduce the results in \cite{btt} where the Hamiltonian case was studied, i.e., 
the case of $G=0$ in (\ref{eq:SGPE}) with $V(x)=|x|^2$ in one dimensional space.  
The authors in \cite{btt} constructed the Gibbs measure, and making use of the invariance of that measure, 
they prove a globalization of the local-in-time solution in a negative Sobolev space, for almost all initial data
with respect to this measure.
We will sometimes make use of their deterministic results. 
Moreover, it was recently proved (see \cite{mtr}) that this Gibbs measure may be obtained as the mean field limit 
of the (finite dimensional) quantum particle system.

\section{Preliminaries and main results}

We consider the following Gross-Pitaevskii equation (complex Ginzburg-Landau equation)  
with a harmonic potential, driven by a space-time white noise in one spatial dimension.   
\begin{equation} \label{eq:CGL}
\left\{
\begin{array}{ll}
dX & = (i+\gamma )(H X +\eta X -\lambda |X|^2 X)dt +\sqrt{2\gamma} dW, \quad t>0, \quad x\in \R,\\
X(0) &= X_0,
\end{array}
\right.
\end{equation}
where $H=\p_x^2-x^2$, $\gamma>0$ and $\eta \ge 0$. 
In this paper we always assume that the nonlinearity is defocusing, namely, we fix $\lambda=1$. 
The unknown function $X$ is a complex valued random field
on a probability space $(\Omega, \mathcal{F}, \prob)$ 
endowed with a standard filtration $(\mathcal{F}_t)_{t\ge 0}$.  
The stochastic process $(W(t))_{t\ge 0}$ is a cylindrical Wiener process on $L^2(\R, \C)$  
associated with the filtration $(\mathcal{F}_t)_{t\ge 0}$, i.e., 
for any complete orthonormal system $\{e_k\}_{k\in \N}$ in $L^2(\R, \R)$, we can write 
$$ W(t,x)=\sum_{k\in \N} \beta_k(t) e_k(x).$$
Here, $\{\beta_k\}_{k\in \N}$ is a sequence of complex-valued independent Brownian motions on the stochastic basis 
$(\Omega, \mathcal{F}, \prob, (\mathcal{F}_t)_{t\ge 0})$. 

It is known that the operator $H$ has a self-adjoint extension on $L^2(\R, \C)$, which is still denoted by $H$, 
and the resolvent of $H$ is compact. Thus, the whole spectrum of $H$ is discrete, and we denote the (real-valued) eigenfunctions by 
$\{h_n\}_{n\ge 0}$ which form an orthonormal basis of $L^2(\R, \R)$. 
They satisfy $H h_n= -\lambda_n^2 h_n$ with $\lambda_n=\sqrt{2n+1}$. In fact, those functions $h_n(x)$ are known as the Hermite functions. 
We denote by $E_N$ the complex vector space spanned by the Hermite functions, $E_N=\mathrm{span}\{h_0, h_1,..., h_N\}$.  
\vspace{3mm}

For $1 \le p \le +\infty,$ and $s\in \R,$ we define the Sobolev space associated to the operator $H$.  
\begin{equation*}
\mathcal{W}^{s,p}(\R)=
\{v \in \mathcal{S}'(\R), \; |v|_{\mathcal{W}^{s,p}(\R)}:=|(-H)^{s/2} v|_{L^p(\R)} <+\infty \},
\end{equation*}
where $\mathcal{S}(\R)$ denotes the Schwartz space. 
If $I$ is an interval of $\R$, $E$ is a Banach space,
and $1\le r\le \infty$, then $L^r(I,E)$ is the space of
strongly Lebesgue measurable functions $v$ from $I$ into $E$
such that the function $t\to |v(t)|_E$ is in $L^r(I)$.
We define similarly the spaces
$C(I,E)$, $C^{\alpha}(I, E)$ or $L^r(\Omega, E)$. 
Sometimes, emphasis will be put in the notation on the measure we consider ; for example if we write 
$L^q ((L^p, d\rho), \R)$, this means the space of real-valued, measurable functions 
on the measure space $(L^p, d\rho)$, with integrable $q$-th power.
For a complex Hilbert space $E$, 
the inner product will be understood as taking the real part, i.e., for $u=u^R+iu^I \in E$ and 
$v=v^R+iv^I \in E$, then we set $(u, v)_{E}:=(u^R, v^R)_{E} + (u^I, v^I)_{E}$    
that is we use the identification $\C \simeq \R^2$. 
\vspace{3mm}

Motivated by the physical background explained in the Introduction, 
our aim is first to show mathematically that Eq.(\ref{eq:CGL}) has the same invariant Gibbs measure as the one described in \cite{btt} 
which is formally described, up to a normalizing constant, by
\begin{equation} \label{formalGibbs}
\rho(du)=e^{-S(u)}du, \quad u\in L^p(\R; \C)
\end{equation}
where $S$ is the Hamiltonian for the case $\gamma=0$, i.e.,
\begin{equation*}
S(u)=\frac{1}{2}\int_{\R} |(-H)^{1/2}u|^2 dx -\frac{\eta}{2} \int_{\R} |u|^2 dx+ \frac{1}{4} \int_{\R}|u|^4 dx,  
\end{equation*}
then to use this in order to prove the global existence of solutions and the exponential convergence to this equilibrium as $t\to +\infty.$
\vspace{3mm}

%

We may give a meaning to the measure (\ref{formalGibbs}) as follows. 
When $u= \sum_{n=0}^{\infty} c_n h_n$, $c_n \in \C$,
we write $c_n= a_n +ib_n$ with $(a_n, b_n) \in \R^2$.
For $N\in \N$, we consider the probability measure on ${\R}^{2(N+1)}$ 
defined by 
\begin{equation*}
d\mu_N := \prod_{n=0}^N \frac{\lambda_n^2}{2\pi} e^{-\frac{\lambda_n^2}{2} (a_n^2+b_n^2)}da_n db_n.    
\end{equation*}
This measure defines a measure on $E_N$ through the map from $\R^{2(N+1)}$ to $E_N$ defined by
\begin{equation*}
(a_n, b_n)_{n=0}^N  \mapsto \sum_{n=0}^N (a_n+ib_n)h_n,
\end{equation*}
and this measure will be again denoted by $\mu_N$. 
Here, $\prod_{n=0}^N \frac{\lambda_n^2}{2\pi}$ is the normalizing factor, i.e., $\mu_N(E_N)=1$.
This measure $\mu_N$ can be seen as the law of the $E_N$-valued random variable 
$\varphi_N(\omega, x)$ defined on $(\Omega, \mathcal{F}, \prob)$ by 
\begin{equation*}
\varphi_N (\omega, x) := \sum_{n=0}^N \frac{\sqrt{2}}{\lambda_n} g_n(\omega)h_n(x), 
\end{equation*}
where $\{g_n (\omega)\}_{n=0}^N$ is a system of independent, complex-valued random variables with the law $\mathcal{N}_{\C}(0,1)$.  
It may be seen, using the asymptotic properties of the $L^p$ norm of the Hermite functions $h_n$, i.e.
\begin{equation*}
|h_n|_{L^p(\R)} \le C_p \lambda_n^{-\frac{1}{6} \theta(p)}, 
\end{equation*}
where 
\begin{equation} \label{def:thetap}
\theta(p)=\left\{
\begin{array}{ll}
1 & \mbox{if}\quad p\ge 4, \\
2-\frac{4}{p} & \mbox{if}\quad 2 \le p \le 4
\end{array}
\right.
\end{equation}
(see Lemma 3.2 of \cite{btt}),   
that $\{\varphi_N\}_{N\ge 0}$ is a Cauchy sequence in $L^2(\Omega, L^p(\R, \C))$ for any $p >2$. 
Thus, the limit $\varphi:=\lim_{N \to \infty} \varphi_N$ is well-defined and 
\begin{equation*}
\varphi(\omega, x) = \sum_{n=0}^{\infty} \frac{\sqrt{2}}{\lambda_n} g_n(\omega)h_n(x).
\end{equation*}

We denote by $\mu$ the measure on $L^p(\R; \C)$, with $p>2$,  induced by this random variable $\varphi(\omega)$, that is, 
for any Borel set $A \subset L^p$,
\begin{equation*}
\mu (A) = \prob(\omega\in \Omega,~ \varphi(\omega) \in A). 
\end{equation*}
Note that the measure $\mu$ can be decomposed into 
\begin{equation*}
\mu =\mu_N \otimes \mu_N^{\perp} 
\end{equation*}
where $\mu_N^{\perp}$ is the law of the random variable on $E_N^{\perp}=\{u \in L^p(\R; \C); \forall v \in E_N, (v, u)_{E_N}=0\}, $
\begin{equation*}
\sum_{n=N+1}^{\infty} \frac{\sqrt{2}}{\lambda_n} g_n(\omega)h_n(x).
\end{equation*}
In the case $\eta=0$, recalling that supp $\mu \subset L^4(\R)\cap L^p(\R)$, for any $p>2$, we see that   \\
$\exp\{-\frac{1}{4}|u|_{L^4(\R)}^4\} \in L^1 (L^p(\R), d\mu)$ for any $p>2$, 
and we can define the Gibbs measure $\rho$ as 
\begin{equation*} 
\rho (du)= {\Gamma}^{-1} \exp \Big\{-\frac{1}{4}|u|_{L^4(\R)}^4 \Big\} \mu (du), \quad \mbox{for} \quad \mu\mbox{-a.e.} \; u
\end{equation*}
where $\Gamma$ is the normalizing constant, that is, $\Gamma=\int_{L^p} e^{-\frac{1}{4}|u|_{L^4}^4} \mu(du)$.
When $\eta$ is positive (and possibly large), we have to use a different decomposition. Consider the spectral projector 
$$\Pi_{N} \Big(\sum_{n=0}^{\infty} c_n h_n \Big):=  \sum_{n=0}^{N} c_n h_n$$
which is bounded in $L^2$ with a bound equal to $1$, and let $N_0=N_0(\eta)$ be such that
$\lambda_{N_0}^2 \le \eta < \lambda_{N_0+1}^2$, with the convention $\lambda_{-1}=0$. 
Then the operator $A:=-H-\eta(I-\Pi_{N_0})$, with $D(A)=D(-H)$,
is clearly a positive, definite, self-adjoint operator, and arguments similar to those above show that
$e^{-\frac12 (Au,u)_{L^2}}du$ defines, up to a normalizing constant, a Gaussian probability measure $\tilde \mu_{\eta}$ on 
$L^p(\R)$ for any $p>2$. Setting then 
$$\tilde{V}_{\eta}(u)=-\frac{\eta}{2}|\Pi_{N_0} u|_{L^2}^2 +\frac{1}{4}|u|_{L^4}^4$$ 
and
$$\tilde{\Gamma}_{\eta}= \int_{L^p}e^{-\tilde{V}_{\eta}(u)} \tilde{\mu}_{\eta}(du),$$
and noting that
there exists a constant $C_{N_0, \eta}>0$ such that $ \frac{1}{4} |u|_{L^4}^4  \ge \tilde{V}_{\eta}(u) \ge \frac{1}{8} |u|_{L^4}^4 -C_{N_0, \eta}$ 
for $\tilde{\mu}_{\eta}$- a.e. $u\in L^p(\R)$, we deduce that $\tilde V_{\eta} \in L^1(L^p, d\tilde \mu_{\eta})$ and we may define the measure $\rho$ as
\begin{equation} \label{gibbs}
\rho(du)
=\tilde{\Gamma}_{\eta}^{-1}e^{-\tilde{V}_{\eta}(u)} \tilde{\mu}_{\eta}(du).
\end{equation}
\vspace{3mm}
Of course, both definitions coincide when $\eta=0$.

Let us now consider the equation (\ref{eq:CGL}). We take $\{h_k, ih_k\}_{k\ge 0}$ 
as a complete orthonormal system in $L^2(\R, \C)$, namely our cylindrical Wiener process is now  
$$ W(t,x)=\sum_{k=0}^{\infty} (\beta_k^R(t)+i\beta_k^I(t)) h_k(x).$$
Here, $(\beta_k^R(t))_{t\ge 0}$ and $(\beta_k^I(t))_{t\ge 0}$ are sequences of real-valued Brownian motions.
\vspace{3mm}

First, we introduce the linear equation
\begin{equation} \label{eq:linear}
dZ = (i+\gamma )H Z dt +\sqrt{2\gamma} dW, \quad t \in \R, \quad x\in \R.
\end{equation}
Note that we consider here that the process $W$ has been extended to the negative time axis.
The stationary solution $Z_{\infty}$ of (\ref{eq:linear}) can be written as
\begin{equation} \label{sol:stationary}
Z_{\infty}(t)= \sqrt{2\gamma} \int_{-\infty}^t e^{(t-s)(i+\gamma)H}dW(s), \quad t\in \R.
\end{equation}
Expanding $W(t)$ as a series, we may write $Z_{\infty}$ as  
\begin{equation*}
Z_{\infty}(t)=\sqrt{2\gamma} \sum_{k=0}^{\infty} \int_{-\infty}^{t} e^{-(t-s)(i+\gamma)\lambda_k^2}
(d\beta_k^R(s)+id\beta_k^I(s) ) h_k(x).  
\end{equation*}
Since $Z_{\infty}$ is stationary,  
$$\sqrt{2\gamma} \int_{-\infty}^{t} e^{-(t-s)(i+\gamma)\lambda_k^2}
(d\beta_k^R(s)+id\beta_k^I(s) )$$ 
has the same law as  
$$\sqrt{2\gamma} \int_{-\infty}^{0} e^{s(i+\gamma)\lambda_k^2}
(d\beta_k^R(s)+id\beta_k^I(s) )$$ 
which is $\mathcal{N}_{\C} \Big(0,\frac{2}{\lambda_k^2}\Big)$. 
Hence we may write 
$$Z_{\infty}(t)=\sum_{k=0}^{\infty} \frac{\sqrt{2}}{\lambda_k} g_k(\omega,t)h_k(x)$$
where $(g_k(\omega,t))_k$ is a family of independent $\mathcal{N}_{\C}(0,1)$, i.e. 
the law $\mathcal{L}(Z_{\infty}(t))$ is equal to the Gaussian measure $\mu.$
\vspace{3mm}

The regularity of $Z_{\infty}(t)$ is given by the following Lemma.
\begin{lem} \label{regularityZ}
Let $T> 0$ be fixed. 
Let $p>2$, $s\in [0,\frac{1}{6} \theta(p))$ and $\alpha \in (0,\frac{1}{12}\theta(p)-\frac{s}{2})$.  
The stationary solution $Z_{\infty}$ of (\ref{eq:linear}) 
has a modification in $C^{\alpha}([0,T],\mathcal{W}^{s,p}(\R))$. Moreover, 
there exists a positive constant $M_{p,T}$ such that
\begin{eqnarray*}
\E\Big(\sup_{t\in [0,T]} |Z_{\infty}(t)|_{L^p} \Big) \le M_{p,T}.
\end{eqnarray*}
\end{lem}

We will give a proof of this lemma in Appendix A. 
Note that this regularity may not be optimal, but is enough for our purpose. 
\vspace{3mm}

Using the regularity of $Z_{\infty}$ a.s. in $C([0,T], L^p(\R))$ for any $p> 2$  
we establish the local existence of the solution to the equation (\ref{eq:CGL}).

\begin{prop} \label{prop:prop1} Let $\gamma>0$, $\eta \ge 0$, $\lambda=1$ and $p \ge 3$. 
Let $T>0.$ Assume $X_0 \in L^p(\R)$. Then there exists a random stopping time $T^*=T^*_{X_0,\omega}>0$, a.s. 
and a unique solution $X(t)$ adapted to $(\mathcal{F}_t)_{t\ge 0}$
of $(\ref{eq:CGL})$ with $X(0)=X_0$, almost surely in
$C([0,T^*), L^p(\R))$. Moreover, we have almost surely, 
$T^*=T$ or 
$\displaystyle{\limsup_{t \to T^{*}}|X(t)|_{L^p(\R)}=+\infty}$. 
\end{prop}
Note that this proposition is also valid for the case $\lambda=-1$, but we focus on the defocusing case $\lambda=1$. 
The Gibbs measure $\rho$ is, in fact, an invariant measure for (\ref{eq:CGL}). 
With the use of this invariant measure, we obtain the global existence of the solution of (\ref{eq:CGL}) for $\rho$-a.e. $X_0$
(or equivalently for $\mu_{\eta}$-a.e. $X_0$) : 

\begin{thm} \label{thm:thm1} Let $\gamma>0$, $\eta \ge 0$, $\lambda=1$ and $p \ge 3$.    
There exists a $\rho$-measurable set $\mathcal{O} \subset L^p(\R)$ such that $\rho(\mathcal{O})=1$, 
and such that for $X_0 \in \mathcal{O}$ 
there exists a unique solution of \eqref{eq:CGL}, $X(\cdot) \in C([0,\infty), L^p(\R))$ a.s. 
\end{thm}

Let $P_t$ be the transition semigroup 	associated with equation (\ref{eq:CGL}) 
which, thanks to Theorem~\ref{thm:thm1},  is well defined and continuous on $L^2((L^p(\R),d\rho), \R)$ 
for any $t \ge 0$ and for $p\ge 3$. We will actually prove that $P_t$ is defined on the set of Borelian bounded
functions on $L^p(\R)$ for any $p \ge 3$ and that it is strong Feller and irreducible. As a consequence, we
will obtain the following theorem.

\begin{thm}\label{t3}
Let $\gamma>0$, $\eta \ge 0$, $\lambda=1$ and $p \ge 3$.    
For any $X_0 \in L^p(\R)$,
there exists a unique solution of \eqref{eq:CGL}, $X(\cdot) \in C([0,\infty), L^p(\R))$ a.s. 
\end{thm}

Finally, we will use the purely dissipative counterpart of equation \eqref{eq:CGL} to prove the exponential
convergence of the transition semi-group, as stated in the following theorem.

\begin{thm} \label{thm:thm2} 
Let $\gamma>0$, $\eta \ge 0$, $\lambda=1$ and $p \ge 3$. 
Let $\phi \in L^2((L^p,d\rho), \R)$, and $\bar{\phi}^{\rho}=\int_{L^p} \phi(y) d\rho(y)$.
Then $u(t, \cdot):=P_t \phi (\cdot)$ converges exponentially to $\bar{\phi}^{\rho}$ in $L^2((L^p,d\rho), \R)$, 
as $t \to \infty$; more precisely,  
$$
\int_{L^p} |u(t,y)-\bar{\phi}^{\rho}|^2 d\rho(y) \le e^{-2\gamma t} \int_{L^p} |\phi(y) -\bar{\phi}^{\rho}|^2 d\rho(y).
$$
\end{thm} 
\vspace{3mm}

\begin{rem}
The statements in the theorems are restricted to the case $p \ge 3$ 
although we only need $p>2$ for the support of the Gaussian measure $\tilde \mu_\eta$.  
This condition comes from  the cubic nonlinearity in (\ref{eq:CGL}); the results are still valid
for more general nonlinear power $|X|^{2\sigma} X$ under the conditions $p \ge 2\sigma+1$ and $p>2$.  
\end{rem}
\vspace{3mm}


In the course of the proof of the above theorems, we will frequently use an approximation by finite dimensional objects.
We thus define here, as in \cite{btt}, for any $p \in [1,\infty]$, a smooth projection operator $S_N : L^p(\R; \C) \to E_N$ by 
\begin{equation*}
S_N \Big(\sum_{n=0}^{\infty} c_n h_n \Big):= \sum_{n=0}^{\infty} \chi \Big(\frac{2n+1}{2N+1}\Big) c_n h_n
=\chi\Big(\frac{H}{2N+1}\Big) \Big(\sum_{n=0}^{\infty} c_n h_n \Big).
\end{equation*}
where $\chi$ is a cut-off function such that $\chi \in C_0^{\infty}(-1,1)$, $\chi=1$ on $[-\frac{1}{2}, \frac{1}{2}].$ 
We will make use of Proposition 4.1 of \cite{btt} : 
$S_N$ is a bounded operator from $L^p$ to $L^p$, uniformly in $N$,  for any $p \in [1,\infty]$. Note that the usual spectral projector 
$\Pi_{N}$ does not satisfy this property.
On the other hand, one can in fact check under the assumptions in Lemma \ref{regularityZ} that  
\begin{equation*}
\E\Big(\sup_{t\in [0,T]} |\Pi_N Z_{\infty}(t)|_{L^p} \Big) \le M_{p,T}, \quad 
\E\Big(\sup_{t\in [0,T]} |(I-\Pi_N) Z_{\infty}(t)|_{L^p} \Big) \le M_{p,T},  
\end{equation*}
and $\Pi_N Z_{\infty}(t)$, $(I-\Pi_N)Z_{\infty}(t)$ make sense in $L^p(\R)$, a.s. if $p>2$ (see the proof of Lemma \ref{regularityZ}).
\vspace{3mm}
 
Let $E$ be a separable Hilbert space and $K$ be a Banach space.   
Given a differentiable function $\varphi$ from $E$ to $K$, we denote by $D\varphi(x)$ its differential at $x\in E.$ 
It is an element of $\mathcal{L}(E,K)$ and if $K=\R$ it is identified with its gradient so that it is also 
seen as an element of $E$. If $\varphi$ is twice differentiable, $D^2 \varphi$ is its second differential. 
Again, we identify $D^2\varphi(x)$, $x\in E$, with an element of $\mathcal{L}(E,E)$ in case of $K=\R$.
If $\{e_i\}_{i\in \N}$ is a Hilbert basis in $E$, $\mu$ a probability measure on $E$, and $\varphi \in C^1_b(E;\R)$, then
$$|D\varphi|_{L^2((E, d\mu), E)}^2:= \int_E |D \varphi(x)|_E^2 d\mu(x)= \int_E \sum_{i\in \N} (D\varphi (x), e_i)_E^2 d\mu(x).$$
\vspace{3mm}

This paper is organized as follows. In Section 3, we will introduce the deterministic properties of the equations 
including the kernel estimate for the deterministic linear part of
equation (\ref{eq:CGL}). The proof of proposition \ref{prop:prop1} will be an easy consequence of
those estimates. A review for the purely dissipative equation, i.e. equation (\ref{eq:CGL}) without the skew-symmetric part 
induced by the imaginary unit $i$ will be given in Section 4 
concerning the existence and the uniqueness of invariant measure, 
and the Poincar\'{e} inequality, for finite dimensional approximations of the equation. The invariance for the finite dimensional
approximations of equation \eqref{eq:CGL} will also be discussed in Section 4.
Note that the results of Section 4 make use of rather standard techniques.
We will establish the global existence of solutions of (\ref{eq:CGL}) for a.e. initial data with respect to 
the Gibbs measure in Section 5, together with the invariance of the (infinite dimensional) measure $\rho$, 
while the strong Feller property and the global existence for all initial data
will be proved in Section 6. Using the purely dissipative equation as an auxiliary equation, 
the convergence to the Gibbs equilibrium in (\ref{eq:CGL}) will be proved in Section 7. Note that, to our knowledge,
it is the first time that such an argument, is used in the infinite dimensional case. In order to simplify the presentation, it will be
assumed from Section 4 to Section 7 that $\eta=0$. Section 8 will be devoted to explain how the arguments can be adapted
to the case $\eta>0$ (and possibly large). In the Appendix, 
we will show some regularity properties of the stochastic convolution needed in the course of the proofs. 
\section{Local existence of the strong solution}

In this section we will give a proof of Proposition \ref{prop:prop1}.
Let $T>0$ and fix $p\ge 3$. Let $z \in C([0,T], L^p)$. 
Consider the following deterministic equation :
\begin{equation} \label{eq:deterministic}
\left\{
\begin{array}{ll}
\partial_t v & = (i+\gamma )(H v +\eta(v+z) -|v+ z|^2 (v+ z)), \quad t>0, \; x\in \R,\\
v(0) &= v_0.
\end{array}
\right.
\end{equation}
Note that $X$ is a solution of \eqref{eq:CGL} if and only if $X=v+z$ with 
$$ z(t) = Z_{\infty}(t)-e^{(i+\gamma)tH} Z_{\infty}(0) =\sqrt{2\gamma} \int_0^t e^{(t-s)(i+\gamma)H}dW(s),$$
and $v$ solution of \eqref{eq:deterministic} with $v_0=X_0$. Moreover, Lemma \ref{regularityZ}, and Lemma \ref{lem:linear_deterministic}
below show that $z$ has paths a.s. in $C([0,T], L^p)$.
Concerning $v$, we can prove the following local existence result.
\begin{prop} \label{prop:prop2}
Let $T>0$ and $p\ge 3$, and let $v_0 \in L^p(\R)$ and $z \in C([0,T], L^p)$. 
Then there exist a time $T^*= T^*_{v_0, z}>0$ 
and a unique solution $v$ 
of (\ref{eq:deterministic})  in $C([0,T^*), L^p(\R))$. 
Moreover, we have 
$T^{*}=T$ or 
$\displaystyle{\limsup_{t \to T^{*}}|v(t)|_{L^p(\R)}=+\infty}$. 
\end{prop}

To prove Proposition \ref{prop:prop2}, we need some estimates for the linear deterministic equation :
\begin{equation} \label{eq:linear_deterministic}
\left\{
\begin{array}{ll}
\p_t w & = (i+\gamma )H w, \quad t>0, \quad x\in \R,\\
w(0) &= f \in \mathcal{S}(\R). 
\end{array}
\right.
\end{equation}

\begin{lem} \label{lem:linear_deterministic}
The solution $w$ of equation (\ref{eq:linear_deterministic}) can be written as 
\begin{equation*}
w(t,x)=e^{t(i+\gamma)H}f= \int_{\R} K_t(x,y) f(y) dy, \quad t>0, \quad x\in \R, 
\end{equation*}
with the kernel 
\begin{equation*}
K_t (x,y):=\frac{1}{\sqrt{-2\pi i \sin(2(\gamma i-1)t)}} 
\exp \Big\{ \frac{\cos(2(\gamma i-1)t)}{i\sin(2(\gamma i -1)t)}\frac{x^2+y^2}{2}
-\frac{1}{i\sin(2(\gamma i -1)t)} xy \Big\}.   
\end{equation*}
The kernel satisfies, for $t>0$ sufficiently small,  
\begin{equation} \label{ineq:young}
|e^{t(i+\gamma)H}f|_{L^r(\R)} \le C_{\gamma} t^{-\frac{1}{2l}}|f|_{L^s(\R)},
\end{equation}
where 
\begin{equation*}
0 \le \frac{1}{r} \le \frac{1}{r} +\frac{1}{l}=\frac{1}{s} \le 1. 
\end{equation*}
\end{lem}
\vspace{3mm}

\begin{rem}
The constant $C_{\gamma}$ in (\ref{ineq:young}) is independent of $t$ and diverges when $\gamma $ is close to $0$. 
\end{rem}

\proof The form of the kernel is due to the Mehler formula (see e.g.\cite{h}). For (\ref{ineq:young}), we decompose the kernel as follows. 
\begin{equation*}
K_t (x,y)= \sqrt{\frac{\delta}{\pi}} e^{-(\beta-\delta)x^2} e^{-\delta (x-y)^2} e^{-(\beta-\delta)y^2},
\end{equation*}
where we put
\begin{eqnarray*}
\delta = -\frac{1}{2} \frac{1}{i\sin (2(\gamma i -1)t)}, \quad 
\beta = -\frac{1}{2} \frac{\cos(2(\gamma i-1)t)}{i\sin (2(\gamma i-1)t)}.
\end{eqnarray*}
Remark that $\mathrm{Re} (\delta) > 0$ for $0< t < \frac{\pi}{4}$ and $\mathrm{Re}(\beta-\delta) > 0$ for any $0< t < \frac{\pi}{4}$;
indeed, we can compute  
\begin{eqnarray*}
i\sin (2(\gamma i -1)t) &=& - \sinh (2\gamma t) \cos (2t) -i \cosh (2\gamma t) \sin (2t), 
\end{eqnarray*}
and thus, $\mathrm{Re}(i\sin (2(\gamma i -1)t))= -\sinh (2\gamma t) \cos (2t) < 0$ for $0 < t < \frac{\pi}{4}$. 
Therefore we get 
$\mathrm{Re} (\delta) > 0$ for $0 < t < \frac{\pi}{4}$ since  
\begin{eqnarray*}
\mathrm{Re} (\delta)= -\frac{1}{2}\mathrm{Re} \Big(\frac{1}{i\sin (2(\gamma i -1)t)} \Big) 
=-\frac12 \frac{\mathrm{Re}(i\sin (2(\gamma i -1)t))}{|\sin (2(i\gamma-1)t)|^2}. 
\end{eqnarray*}
On the other hand, 
\begin{eqnarray*}
\beta-\delta =-\frac{1}{2} \frac{\cos(2(\gamma i-1)t)-1}{i\sin (2(\gamma i-1)t)}. 
\end{eqnarray*}
Since 
\begin{eqnarray*}
\cos (2(\gamma i -1)t) &=& \cosh (2\gamma t) \cos (2t) +i \sinh (2\gamma t) \sin (2t), 
\end{eqnarray*}
we have, 
\begin{eqnarray*}
&& \frac{\cos(2(\gamma i-1)t)-1}{i\sin (2(\gamma i-1)t)}\\
&& =-\frac{(\cosh(2\gamma t) \cos (2t)-1+i\sinh (2\gamma t) \sin (2t))(\sinh (2\gamma t) \cos(2t)-i\cosh(2\gamma t)\sin (2t))}
{\sinh^2 (2\gamma t) \cos^2 (2t)+\cosh^2 (2\gamma t) \sin^2 (2t)}.
\end{eqnarray*}
The real part of the numerator is 
\begin{eqnarray*}
&& -(\cosh(2\gamma t) \cos (2t)-1)(\sinh (2\gamma t) \cos(2t))-\sinh (2\gamma t) \cosh(2\gamma t)\sin^2 (2t)\\
&&= -\sinh(2\gamma t)[-\cos (2t)+\cosh (2\gamma t)] \\ 
&&< 0,
\end{eqnarray*}
for $0< t < \frac{\pi}{4}$ since $\cosh (2\gamma t) \ge 1 > \cos(2t)$. Accordingly, $\mathrm{Re}(\beta-\delta)> 0$ for any $0< t < \frac{\pi}{4}$.
\vspace{3mm}

Then, we have,
\begin{eqnarray*}
|e^{t(i+\gamma)H} f|_{L^r(\R)} &\le &\Big\{ \int_{\R} \Big(\int_{\R} |K_t(x,y) f(y)| dy \Big)^r dx \Big\}^{1/r}\\
&\le &\Big\{ \int_{\R} \Big(\int_{\R} 
\sqrt{\frac{|\delta|}{\pi}} e^{-\mathrm{Re}(\beta-\delta)x^2} e^{-\mathrm{Re}\delta (x-y)^2} 
e^{-\mathrm{Re}(\beta-\delta)y^2} |f(y)| dy \Big)^r dx \Big\}^{1/r} \\
&\le & \sqrt{\frac{|\delta|}{\pi}} \Big|e^{-(\mathrm{Re}\delta) x^2} \ast |f| \Big|_{L^r(\R)} \\
&\le & \sqrt{\frac{|\delta|}{\pi}} |e^{-(\mathrm{Re}\delta) x^2}|_{L^{l'}(\R)} |f|_{L^s(\R)}, 
\quad \mbox{with} \quad \frac{1}{r} +1=\frac{1}{l'}+\frac{1}{s},
\end{eqnarray*}
where we have used the Young inequality in the last line. Then, it suffices to show that for $0 < t< \frac{\pi}{4}$,
\begin{equation}
\label{dec}
 \sqrt{\frac{|\delta|}{\pi}} |e^{-(\mathrm{Re}\delta) x^2}|_{L^{l'}(\R)} \le C t^{-\frac{1}{2l}}
\end{equation}
for $l\ge 1$ such that $\frac{1}{l'}+\frac{1}{l}=1$. 
Note first that there exists a constant $C>0$ which is independent 
of $t \in (0,\frac{\pi}{4})$ such that 
\begin{equation*}
\Big|\sqrt{\frac{|\delta|}{\pi}} e^{-(\mathrm{Re}\delta) x^2}\Big|_{L^{\infty}(\R)} \le \sqrt{\frac{|\delta|}{\pi}} \le C t^{-1/2}, 
\end{equation*}
for $0 < t < \frac{\pi}{4}.$
Next, note that for sufficiently small $t>0$,
$$
\Big| \frac{\mathrm{Im}\delta}{\mathrm{Re}\delta}\Big| = \frac{\sin 2t}{\sinh (2\gamma t)}  \frac{\cosh(2\gamma t)}{\cos2t}
\le \frac{(e^{4\gamma^2 t}+1)(\sin 2t)}{16 \gamma^4 t \cos(2t)}
\le \frac{1}{2\gamma^4}
$$
which is a positive constant if $\gamma \ne 0$, 
and we deduce from this that for some constant $C_{\gamma}$,
$$
\Big|\sqrt{\frac{|\delta|}{\pi}} e^{-(\mathrm{Re}\delta) x^2}\Big|_{L^{1}(\R)}=\sqrt{\frac{|\delta|}{\mathrm{Re}\delta}}\le C_{\gamma}.
$$
Inequality \eqref{dec} follows by interpolation between the cases $l=1$ and $l=+\infty$.
\hfill \qed  

\proof (of Proposition \ref{prop:prop2}). We follow the arguments in \cite{gv}. 
Let $T_0 \le T$ be small enough for the inequality (\ref{ineq:young}) in Lemma \ref{lem:linear_deterministic} to be satisfied.   
We consider the closed ball in $C([0,T_0], L^p)$,
\begin{equation*}
B_R(T_0):= \{v \in C([0,T_0], L^p),~ |v|_{C([0,T_0],L^p)} \le R\} 
\end{equation*}
with $R:=2C_{\gamma}|v_0|_{L^p}$ and $C_{\gamma}$ is the constant appearing in (\ref{ineq:young}). 
We will check that the map $\mathcal{T}$ defined by 
\begin{equation*}
\mathcal{T}v (t):= e^{t(i+\gamma)H}v_0 - (i+\gamma) \int_0^t e^{(t-s)(i+\gamma)H}\left[|v+z|^2-\eta\right] (v+z)(s) ds 
\end{equation*}
is a strict contraction on $B_R(T_0)$, for a possibly smaller $T_0$.
For any $v_1, v_2 \in B_R(T_0),$ taking $r=p$ and $s=p/3$ in (\ref{ineq:young}), 
\begin{eqnarray*}
|\mathcal{T}(v_1)-\mathcal{T}(v_2)|_{L^p} 
&\le& C_{\gamma} \int_0^t |t-s|^{-\frac{1}{p}}|\left[|v_1+z|^2 (v_1+z)-|v_2+z|^2 (v_2+z)\right](s)|_{L^{p/3}} ds\\
& & \hskip 0.1 in + \; C_{\gamma} \eta \int_0^t |v_1(s)-v_2(s)|_{L^p} ds\\
&\le & C_{\gamma} \int_0^t \left[|t-s|^{-\frac{1}{p}} \max_{i=1,2}|v_i(s)+z|_{L^p}^2 +\, \etaÊ\right] |v_1(s)-v_2(s)|_{L^{p}}ds,
\end{eqnarray*}
where we have used H\"{o}lder inequality in the last inequality. 
Put $\theta :=1-\frac{1}{p}>0$. Then
\begin{eqnarray*}
|\mathcal{T}(v_1)-\mathcal{T}(v_2)|_{C([0,T_0],L^p)} 
&\le& C_{\gamma} \left[ T_0^{\theta}  
\max_{i=1,2}|v_i+z|_{C([0,T_0],L^{p})}^2 +\, \eta T_0\right] |v_1-v_2|_{C([0,T_0],L^{p})}\\
&\le& C_{\gamma} \left[T_0^{\theta} (R^2+|z|_{C([0,T],L^{p})}^2) +\, \eta T_0 \right] |v_1-v_2|_{C([0,T_0],L^{p})}. 
\end{eqnarray*}
Using now the inequality (\ref{ineq:young}) with $r=s=p$ for the free term, we similarly have, 
for $v_1 \in B_R(T_0),$
\begin{eqnarray*} 
|\mathcal{T}(v_1)|_{C([0,T_0],L^p)} 
&\le& C_{\gamma} |v_0|_{L^p} +C_{\gamma} T_0^{\theta}(R^3+|z|_{C([0,T_0],L^{p})}^3)+\, \eta T_0 \\
&\le& \frac{R}{2} +C_{\gamma} T_0^{\theta}(R^3+|z|_{C([0,T],L^{p})}^3)+\, \eta T_0.
\end{eqnarray*}
We thus see that $\mathcal{T}$ is a contraction in $B_R(T_0)$ provided $T_0$ is sufficiently small, 
which gives a unique solution in $C([0,T_0], L^p)$, and, thanks to classical extension arguments,
a unique maximal solution in $C([0,T^*), L^p)$, where $T^*$ depends only on $\gamma,$ $\eta,$ $v_0$ and $z$. 
\hfill\qed
\vspace{3mm}

Next, we introduce an approximation for the solution $v$ of (\ref{eq:deterministic}) defined above. 
We fix $T>0$, and $p\ge 3$. 
For any $z \in C([0,T], L^p(\R)\cap L^4(\R))$, let $v_N$ be the solution of  
\begin{equation} \label{eq:deterministic_Galerkin1}
\p_t v = (i+\gamma )\big(H v +\eta \, S_N(v+z) - S_N ( |S_N (v + z)|^2 S_N(v+ z))\big), \quad t>0, \quad x\in \R, 
\end{equation}
with initial data
\begin{equation}
\label{eq:deterministic_Galerkin2}
v_N(0)\in E_N. 
\end{equation}

This equation has an energy estimate, which is easily obtained by taking the $L^2$- inner product of \eqref{eq:deterministic_Galerkin1}
with $v_N$ and using the boundedness of $S_N$ in $L^4(\R)$ and Young's inequality (see also \cite{bs,gv}):
\begin{eqnarray} \label{ineq:energy} 
& & \frac{1}{2} \frac{d}{dt} |v_N(t)|_{L^2}^2 +\gamma |v_N(t)|_{\mathcal{W}^{1,2}(\R)}^2 - \eta |v_N(t)|_{L^2(\R)}^2
+\frac{\gamma}{2} \int_{\R} |v_N(t,x)|^4 dx \\ \nonumber & & \le C_{\gamma} \int_{\R} |S_N z (t)|^4 dx \\ \nonumber & &
\le C_{\gamma}\sup_{t\in [0,T]} |z (t)|_{L^4}^4 \le C_{\gamma,T}.
\end{eqnarray}
Therefore, $|v_N(t)|_{L^2} \le e^{2(\eta-\lambda_0^2 \gamma)T}\left[|v_N(0)|_{L^2} +C_{\gamma, T, z}\right]$ for all $t \in [0,T]$. 
Since $v_N(0) \in E_N$, any norm of $v_N(0)$ is finite, i.e. $|v_N(0)|_{L^r} \le C_{N,r}$ for any $r\in [1,\infty]$, 
and for any $p\ge 1,$
$$|v_N(t)|_{L^p} \le C_{N,p}|v_N(t)|_{L^2} \le C_{N,p}e^{2(\eta-\lambda_0^2 \gamma)T}\left(|v_N(0)|_{L^2} +C_{\gamma, T, z}\right)\le C_{p,N,\gamma,T,z}$$ 
for all $t \in [0,T]$. 
\vspace{3mm}

Remark also that by Proposition 4.1 of \cite{btt}, the bound of $S_N$ from $L^p$ to $L^p$ is uniform in $N$, 
thus for $z \in C([0,T], L^p(\R))$, 
\begin{equation} \label{eq:convergence_S_N}
\lim_{N\to \infty}|(S_N-I) z|_{L^{\infty}(0,T, L^p(\R))}=0, \quad  1 \le p \le +\infty.
\end{equation}
We thus have the following proposition. 

\begin{prop} \label{prop:prop4} 
Let $p\ge 3$ and $T>0$. Let $z\in C([0,T], L^p(\R) \cap L^4(\R)).$ Then there exists 
a unique global solution in $C([0,T], E_N)$, denoted by $v_N$, of (\ref{eq:deterministic_Galerkin1})-(\ref{eq:deterministic_Galerkin2}). 
Moreover, for $v_0 \in L^p(\R)$. If $v_N(0)=S_N v_0$, then
$v_N$ converges to $v$ in $L^{\infty}(0,T, L^p(\R))$ as $N \to \infty$ for any $T<T^*$.
\end{prop}

We remark that the convergence of $v_N$ to $v$ is a consequence of the fixed point argument seen in the proof of Proposition 
\ref{prop:prop2} and of (\ref{eq:convergence_S_N}). 

\section{Purely dissipative case and Gibbs measure}

We recall in this section a few results concerning the purely dissipative case that will be useful in the following.
Although those results are obtained thanks to fairly standard techniques, we will give some details for the sake of completeness.
As explained in  Section 2, we assume from now on that $\eta=0$, in order that the arguments are simpler to state, and we refer to
Section 8 for a description of the case $\eta>0$.
\vspace{3mm}

The purely dissipative equation is then
\begin{equation} \label{equ:diss}
dY=\gamma (HY-|Y|^2 Y)dt+\sqrt{2\gamma} dW. 
\end{equation}
We first consider the following approximation :
\begin{equation} \label{equ:dissipative}
dY =\gamma (HY-S_N (|S_N Y|^2 S_N Y))dt + \sqrt{2\gamma} \Pi_N dW, \quad Y(0)=y \in E_N.
\end{equation}
If we put $u_N:=Y_N-\Pi_N Z_{\infty}$ where $Y_N$ satisfies (\ref{equ:dissipative}), then $u_N$ verifies
\begin{equation} \label{eqforu}
\partial_t u  = \gamma \left(H u - S_N(|S_N(u+ Z_{\infty})|^2 S_N(u+ Z_{\infty}))\right),  \; x\in \R.
\end{equation}
Recall that $Z_{\infty}$ is defined in (\ref{sol:stationary}), and note that $S_N \circ \Pi_N =S_N$.  
By the same arguments as for (\ref{eq:deterministic_Galerkin1}), replacing $z$ by $\Pi_N Z_{\infty}$, 
 and using Lemma \ref{regularityZ}, we have the energy estimate (\ref{ineq:energy}) for $u_N$, almost surely.
It thus follows as in Proposition \ref{prop:prop4} that for any $\gamma \ne 0$ and $p\ge 3$, there exists 
a unique global solution $u_N \in C(\R^+, E_N)$ a.s. of equation \eqref{eqforu} so that 
$Y_N=u_N+\Pi_N Z_{\infty} \in C(\R^+, E_N)$, a.s.
We denote by $Y_N(t,0,y)$ the solution of \eqref{equ:dissipative} with initial data $Y_N(0,0,y)=y$. 
In fact, since the noise has been extended to the negative time axis, by the same arguments as above, 
we may consider $Y_N(t, -t_1, y)$ with a $t_1>0$.  

The dissipativity inequality, i.e. the inequality
\begin{eqnarray*}
(\gamma H(y-z) -\gamma S_N(|S_N y|^2 S_N y -|S_N z|^2 S_N z), y-z)_{E_N} \le -\gamma |y-z|_{E_N}^2, 
\end{eqnarray*}
which holds for all $y,z \in E_N$, and for all $\gamma>0$, implies that for any $y,z \in E_N$, and for any $t \ge -t_1$,  
\begin{equation} \label{ineq:dissipativity} 
|Y_N(t,-t_1,y)-Y_N(t,-t_1,z)|_{E_N} \le e^{-\gamma (t+t_1)} |y-z|_{E_N}.
\end{equation}
Using \eqref{ineq:dissipativity} with $z=0$ together with Lemma \ref{regularityZ} (whose proof is valid also for negative time intervals),
one easily obtains for any $y \in E_N$, and for any $t \ge -t_1$, 
$$
\E(|Y_N(t,-t_1,y)|_{E_N}) \le C + |y|_{E_N}.
$$
We deduce for all $t \ge -t_1 \ge -t_2$ with $t_1, t_2>0$,
\begin{eqnarray*}
\E (|Y_N(t,-t_1,y)-Y_N(t,-t_2,y)|_{E_N}) &\le & e^{-\gamma (t+t_1)} \E(|Y_N(-t_1,-t_2,y)-y|_{E_N})\\
& \le & e^{-\gamma(t+t_1)} (2|y|_{E_N}+C). 
\end{eqnarray*}
It follows, in particular letting $t=0$, that the family $\{Y_N(0,s,y), s\le 0\}$ is a Cauchy family  in $L^1(\Omega, E_N)$ as $s \to -\infty$,
and there exists a random variable $\xi_N \in L^1(\Omega, E_N)$ 
such that $\xi_N:=\lim_{s\to -\infty}Y_N(0,s,y)$ in 
$L^1(\Omega, E_N)$. Thus, as $t \to +\infty$, 
$$ \mathcal{L}(Y_N(0,-t,y)) \to \mathcal{L}(\xi_N),$$
in the weak topology. 
Then, setting $\bar{\rho}_N:=\mathcal{L}(\xi_N),$ it follows that $\bar{\rho}_N$ is the unique invariant measure of the flow
generated by (\ref{equ:dissipative}). 
\vspace{3mm}

More details for the purely dissipative case are found in, e.g., \cite{d,dpz0,ms}.  
\vspace{3mm}

\begin{rem}
In the purely dissipative case, the solution $Y$ of \eqref{equ:diss} can directly be shown to be global, since  
one may prove a $L^p$-bound on $u=Y-Z_{\infty}$ for any $\gamma> 0$, $p\ge 4$ 
(see \cite{bs,gv}) :
\begin{equation} \label{est:Lp}
|u(t)|_{L^p(\R)} \le C_T, \quad t\in [0,T]
\end{equation}
for any $T>0.$
\end{rem}
\vspace{3mm}

Now let our attention turn to the Gibbs measure.
Define for $y\in E_N$,  
\begin{equation} \label{finiteGibbs}
d\tilde{\rho}_N(y):=\tilde{\Gamma}_N^{-1} e^{-\frac{1}{4}|S_N y|_{L^4}^4} d \mu_N(y), 
\quad \tilde{\Gamma}_N = \int_{E_N} e^{-\frac{1}{4}|S_N y|_{L^4}^4} d\mu_N(y). 
\end{equation}

The next proposition will state that $\tilde{\rho}_N$ is also an invariant measure 
for the flow generated by (\ref{equ:dissipative}), which implies $\bar{\rho}_N=\tilde{\rho}_N$ by the above uniqueness
property.
On the other hand, we will see that the measure $\tilde{\rho}_N$ 
is also  invariant for the flow of the approximation to (\ref{eq:CGL}) :
\begin{equation} \label{eq:CGL_finite_Galerkin}
dX =(i+\gamma)(HX-S_N (|S_N X|^2 S_N X))dt + \sqrt{2\gamma} \Pi_N dW, \quad X(0) \in E_N. 
\end{equation}
Remind that the solution of Equation (\ref{eq:CGL_finite_Galerkin}) exists globally in $C(\R^+, E_N)$ a.s.;
such a solution will be denoted by $X_N$. Indeed, $X_N$  
can be written as $X_N=v_N + \Pi_N Z_{\infty}$ where $v_N$ is the solution of (\ref{eq:deterministic_Galerkin1}) 
with $z=\Pi_N Z_{\infty}$ and $v_N (0)=X_N (0)-\Pi_N Z_{\infty}(0) \in E_N$ in (\ref{eq:deterministic_Galerkin2}). 

\begin{prop} \label{invariance}
The measure $\tilde{\rho}_N$ defined in (\ref{finiteGibbs}) is invariant 
by the flow of (\ref{equ:dissipative}) and by the flow of (\ref{eq:CGL_finite_Galerkin}). 
\end{prop}

\proof Write Equation (\ref{eq:CGL_finite_Galerkin}) in the form:
\begin{equation*}
dX = -J DI(X)dt-\gamma DI(X)dt + \sqrt{2\gamma} \Pi_N dW, \quad X(0)=y \in E_N
\end{equation*}
where 
\begin{equation*}
J=
\begin{pmatrix} 
0 & -1 \\
1 & 0 
\end{pmatrix}
:\R^2 \to \R^2, 
\end{equation*}
and 
\begin{equation*}
I(y)= \frac{1}{2} \int_{\R} |(-H)^{1/2} y|^2 dx + \frac{1}{4} \int_{\R}|S_N y|^4 dx, \quad y\in E_N.  
\end{equation*}

This proof will be valid also in the purely dissipative case (\ref{equ:dissipative}),
removing the term $J DI(X)dt$.   
The generator of the transition semigroup for the flow of (\ref{eq:CGL_finite_Galerkin}) is given by
\begin{equation*}
(\mathcal{L}_N f)(y)=\gamma \mathrm{Tr} D^2 f(y)-\gamma (D f(y), DI(y))_{E_N} -(D f(y), J DI(y))_{E_N}
\end{equation*}
for any $f\in C^2_b(E_N, \R^2)$. Here, if $f=(f_1, f_2),$ then
$$\mathrm{Tr} D^2 f(y) = \sum_{k=0}^N \left[(D^2 f_1(y)h_k, h_k)_{E_N}+(D^2 f_2(y) h_k, h_k)_{E_N}\right].$$ 
We shall show that $\tilde{\rho}_N$ is invariant for (\ref{eq:CGL_finite_Galerkin}), that is,  for any $f\in C^2_b(E_N, \R^2)$,
$$ \int_{E_N} (\mathcal{L}_N f) (y) \tilde{\rho}_N (dy)=\int_{E_N} f(y) \mathcal{L}^*_N  \tilde{\rho}_N (dy)=0.$$
Indeed, by direct calculations, 
\begin{eqnarray*}
& &\tilde \Gamma_N \int_{E_N} (\mathcal{L}_N f) (y) \tilde{\rho}_N (dy) \\
&=& \int_{E_N} \Big[\gamma \mathrm{Tr} D^2 f(y)-\gamma (D f(y), DI(y))_{E_N} 
 -(D f(y), J DI(y))_{E_N} \Big] e^{-I(y)} dy \\
&=& \gamma \int_{E_N} \mathrm{Tr} D^2 f(y) e^{-I(y)} dy
+\gamma \int_{E_N} \Big(D f(y), D (e^{-I(y)})\Big)_{E_N} dy \\
 & & - \int_{E_N} f(y)  \mathrm{Tr}(DJD)(e^{-I(y)}) dy\\
&=& 0,
\end{eqnarray*}
where we have used integrations by parts, and the fact that $\mathrm{Tr}(DJD)=0$.
\hfill\qed

\begin{prop} \label{prop:poincare} 
For any $\phi \in C^1_b(E_N)$ 
the following inequality is satisfied :
\begin{eqnarray} \label{ineq:poincare}
\int_{E_N} |D\phi(y)|_{E_N}^2 d\tilde{\rho}_N(y) \ge  \int_{E_N} |\phi(y)-\bar{\phi}^{\tilde{\rho}_N}|^2 d\tilde{\rho}_N(y), 
\end{eqnarray}
where $\bar{\phi}^{\tilde{\rho}_N}=\int_{E_N} \phi(y) d\tilde{\rho}_N(y).$
\end{prop}

\begin{rem} 
The important point here is the fact that the properties of the Gibbs measure $\tilde{\rho}_N$,
which is the invariant measure for (\ref{eq:CGL_finite_Galerkin}) 
may be investigated through the equation (\ref{equ:dissipative}). 
\end{rem}

\proof (of Proposition \ref{prop:poincare}.) Note that we may assume, replacing if necessary $\phi$ by $\phi-\bar\phi^{\tilde{\rho}_N}$, that $\bar\phi^{\tilde{\rho}_N}=0$.
Let $Y_N$ be the unique global solution of (\ref{equ:dissipative}), and 
let $R^N_t$ be the transition semigroup on $E_N$ for (\ref{equ:dissipative}), i.e., 
$R^N_t \phi (y) := \E (\phi(Y_N(t,y)))$ for any 
$\phi \in C_b(E_N, \R)$. Recall that $y=Y_N(0)\in E_N$. Assume now that $\phi \in C^1_b(E_N, \R)$. Then, 
for $h\in E_N$ and $y\in E_N$, 
\begin{equation*}
(D(R^N_t \phi)(y), h)_{E_N}=\E(D\phi(Y_N(t,y)), \eta^h_N)_{E_N},
\end{equation*}
where $\eta^h_N := D_y Y_N(t,y) . h \in E_N$ satisfies 
\begin{equation*} 
\left\{
\begin{array}{ll}
\frac{d\eta}{dt} &= \gamma H \eta -\gamma S_N(|S_N Y_N|^2 S_N \eta) -2\gamma S_N (\mathrm{Re}( \overline{S_N Y_N}  S_N \eta) S_N Y_N), \\ 
\eta(0) &= h \in E_N.
\end{array}
\right.
\end{equation*}
Taking the scalar product in $E_N$ with $\eta^h_N$, we have 
\begin{eqnarray*}
\frac{1}{2} \frac{d}{dt} |\eta^h_N|_{E_N}^2 + \gamma \lambda_0^2 |\eta^h_N|^2_{E_N} 
&\le& -\gamma (|S_N Y_N|^2 S_N \eta^h_N, S_N \eta^h_N)_{E_N} 
-2 \gamma \int_{\R} \Big\{\mathrm{Re} (\overline{S_N Y_N} S_N \eta) \Big\}^2 dx \\
&\le& 0.
\end{eqnarray*}
Thus
\begin{equation}
\label{estimeta}
|\eta_N^h(t)|_{E_N}^2 \le e^{-2\omega t}|h|_{E_N}^2, \quad  \omega = \gamma \lambda_0^2=\gamma,
\end{equation}
which leads to 
\begin{eqnarray*}
|(DR_t^N \phi (y), h)_{E_N}| &\le& \E \Big( |(D\phi(Y_N(t,y)), \eta^h_N)_{E_N}|\Big)\\ 
  &\le& \E (|D\phi(Y_N(t,y))|_{E_N} |\eta_N^h|_{E_N}) \\
  &\le& e^{-\gamma t} |h|_{E_N} \E (|D\phi(Y_N(t,y))|_{E_N}).
\end{eqnarray*}
Put $g(t,y):=R^N_t \phi(y)=\E(\phi(Y_N(t,y)))$, and $\psi(y):=|D\phi(y)|_{E_N}$ for any $y\in E_N$. 
Then $\psi \in C_b(E_N, \R)$ and $(D_yg)(t,y)=D(R^N_t\phi)(y)$. By the above computation, 
$|Dg(y)|_{E_N} \le e^{-\gamma t}|R^N_t \psi(y)|$. Therefore,  we obtain
\begin{eqnarray} \nonumber
\int_{E_N} |D_yg(t,y)|_{E_N}^2 d\tilde{\rho}_N(y) 
&\le& e^{-2\gamma t}\int_{E_N} |R^N_t \psi(y)|^2 d\tilde{\rho}_N(y) \\ \nonumber
&\le& e^{-2\gamma t}\int_{E_N} R^N_t \psi^2 (y) d\tilde{\rho}_N(y) \\ \nonumber
&=& e^{-2\gamma t}\int_{E_N} \psi^2 (y) d\tilde{\rho}_N(y) \\ \label{Dphi}
&=&  e^{-2\gamma t}\int_{E_N} |D\phi(y)|_{E_N}^2 d\tilde{\rho}_N(y), 
\end{eqnarray}
where we have used Cauchy-Schwarz inequality, 
and the invariance of the measure $\tilde{\rho}_N$ for the transition semigroup $R^N_t$.  
On the other hand,  
\begin{equation} \label{eq:dissipative1}
\frac{d}{dt}|g(t,y)|_{L^2(E_N, d\tilde{\rho}_N)}^2= -2\gamma \int_{E_N} |Dg(y)|_{E_N}^2 d\tilde{\rho}_N (y);   
\end{equation}
indeed, let $\mathcal{M}_N$ be the generator of $R^N_t$, so that
$\frac{d}{dt} g=\mathcal{M}_N g$ ; 
then for $y\in E_N,$  by the Kolmogorov equation,
\begin{equation*}
(\mathcal{M}_Ng)(y)=\gamma (\mathrm{Tr} D^2 g)(y)+(D g(y), \gamma (Hy-S_N (|S_N y|^2 S_N y))_{E_N}. 
\end{equation*}
Therefore, 
\begin{equation*}
\mathcal{M}_N g^2 (y)=2\gamma |D g(y)|_{E_N}^2 +2g(y) (\mathcal{M}_Ng)(y). 
\end{equation*}
Then, we have 
\begin{eqnarray*}
0 &=& \int_{E_N} g^2(y) \mathcal{M}_N^*\tilde{\rho}_N(dy)=  \int_{E_N}   \mathcal{M}_N g^2(y) \tilde{\rho}_N(dy) \\
  &=& 2\gamma \int_{E_N} |D g(y)|_{E_N}^2 d \tilde{\rho}_N(y) +2 \int_{E_N} g(y) (\mathcal{M}_Ng)(y) d \tilde{\rho}_N(y)\\
  &=& 2\gamma \int_{E_N} |D g(y)|_{E_N}^2 d \tilde{\rho}_N(y) + \frac{d}{dt} \int_{E_N} |g(t)|^2 d \tilde{\rho}_N(y), 
\end{eqnarray*}
so that \eqref{eq:dissipative1} holds.
Integrating \eqref{eq:dissipative1}  on $[0,t]$, and using (\ref{Dphi}) we have 
\begin{equation} \label{eq:dissipative2}
|g(t,\cdot)|_{L^2(E_N, d\tilde{\rho}_N)}^2-|g(0)|_{L^2(E_N, d\tilde{\rho}_N)}^2 \ge 
-(1-e^{-2\gamma t}) \int_{E_N} |D \phi(y)|_{E_N}^2 d\tilde{\rho}_N (y).
\end{equation}
Note that $|g(0)|_{L^2(E_N, d\tilde{\rho}_N)}^2=|\phi|_{L^2(E_N, d\tilde{\rho}_N)}^2$.
We conclude the proof of (\ref{ineq:poincare}) by letting $t$ tend to  $+\infty$ in (\ref{eq:dissipative2}),
using the following Lemma, whose proof easily follows from the dissipative inequality (\ref{ineq:dissipativity}).
\hfill\qed

\begin{lem} \label{lem:cvgence}
Let $\phi \in C^1_b(E_N, \R)$, then there is a constant $C_N$ such that
$$
|R^N_t \phi -\bar\phi^{\tilde{\rho}_N}|_{L^2(E_N, d\tilde{\rho}_N)} \leq C_N e^{-\gamma t} |\phi|_{\mathrm{Lip}}.
$$
\end{lem}

\proof
Since $\tilde{\rho}_N$ is invariant for $R^N_t$ and $\tilde{\rho}_N(E_N)=1$, we have  for any $y\in E_N$,
\begin{eqnarray*}
 |\E(\phi(Y_N(t,y)))-\bar{\phi}^{\tilde{\rho}_N}|
&=& |\E(\phi(Y_N(t,y)))-\int_{E_N} \phi(z) d\tilde{\rho}_N(z)| \\
&=&|\E(\phi(Y_N(t,y)))-\int_{E_N} R^N_t \phi(z) d\tilde{\rho}_N(z)| \\
&=&|\E \int_{E_N}(\phi(Y_N(t,y))-\phi(Y_N(t,z))) d\tilde{\rho}_N(z) | \\
&\le& e^{-\gamma t} |\phi|_{\mathrm{Lip}} \int_{E_N} |y-z|_{E_N} d\tilde{\rho}_N(z),
\end{eqnarray*}
by (\ref{ineq:dissipativity}) with $t_1=0$. 
Thus,  using Cauchy-Schwarz inequality,
\begin{eqnarray*}
|R^N_t \phi-\bar{\phi}^{\tilde{\rho}_N}|_{L^2(E_N, d\tilde{\rho}_N)}^2
&\le& |\phi|_{\mathrm{Lip}}^2 e^{-2\gamma t} \int_{E_N} \int_{E_N} |y-z|_{E_N}^2 d\tilde{\rho}_N(y)d\tilde{\rho}_N(z) \\
&\le& 2 |\phi|_{\mathrm{Lip}}^2 e^{-2\gamma t} \int_{E_N} \int_{E_N} (|y|^2_{E_N} + |z|^2_{E_N}) d\tilde{\rho}_N(y)d\tilde{\rho}_N(z).
\end{eqnarray*}
The integral on the right hand side is finite, since $\mu_N$ is a Gaussian measure on $E_N$, 
$d\tilde \rho_N= \Gamma_N^{-1} \exp(-|S_N y|_{L^4}^4/4) d\mu_N$
and $\exp(-|S_N y|_{L^4}^4/4) \le 1$. The result follows.
\hfill\qed

\section{Almost sure global existence}

A similar idea to \cite{dpd} will be used to prove Theorem \ref{thm:thm1}.  
Let $X_0 \in L^p$ with $p\ge 3$ be given.  
Consider a particular decomposition for the solution $X$ of (\ref{eq:CGL}), given by $X(t)=v(t)+Z_{\infty}(t)-Z_{\infty}(0)$, 
where $v$ is the solution of (\ref{eq:deterministic}) with $z(t)=Z_{\infty}(t)-Z_{\infty}(0)$, $v_0=X_0$. 
Recall that $Z_{\infty}$ is the stationary solution for (\ref{eq:linear}), which is written as (\ref{sol:stationary}).  
Recall also that the law of $Z_{\infty}(t)$, supported on $L^p$ if $p>2$, is equal to $\mu$.  
We will see that the measure $\rho$ in (\ref{gibbs}) should be invariant for the flow of $X$,  
as the limit of 
\begin{eqnarray}
\label{defrhoN}
d\rho_N(u)&=& {\Gamma}_N^{-1} \exp\Big\{-\frac{1}{4}|S_N u|_{L^4}^4\Big\} 
d\mu_N(\Pi_N u)d\mu_N^{\perp}((I-\Pi_N)u), \quad u\in L^p(\R),
\end{eqnarray}
where $p\ge 3$ and 
\begin{equation*}
\Gamma_N=\int_{L^p} \exp\Big\{-\frac{1}{4}|S_N u|_{L^4}^4\Big\} d\mu(u).
\end{equation*}
This probability measure $\rho_N$ is an invariant measure for the flow of a certain approximation $X^N$ to $X$. Making use of this invariance,  
we will show that $X$ exists globally in $L^p(\R)$ for $p\ge 3$. 
\vspace{3mm}

Let us first consider the approximation $X_N$ given by the solution of (\ref{eq:CGL_finite_Galerkin}), 
which may be written as $X_N(t)=v_N (t)+ \Pi_N Z_{\infty}(t)-\Pi_N Z_{\infty}(0)$. Here, $v_N$ is the solution 
of (\ref{eq:deterministic_Galerkin1}) choosing $v_N(0)=S_N X_0.$ Namely, $X_N(0)=S_N X_0.$ 
Now, we define $X^N(t)=X_N(t)+\tilde{Z}_N(t),$ where $\tilde{Z}_N$ satisfies
\begin{equation*}
\left\{
\begin{array}{ll}
dZ = (i+\gamma)HZ+\sqrt{2\gamma}(I-\Pi_N)dW, &\\
Z(0)=(I-S_N) X_0.&
\end{array}
\right.
\end{equation*}
Then, noting that $S_N \tilde Z_N =0$, one sees that $X^N$ satisfies the following equation with initial data $X^N(0)=X_0$ :
\begin{equation} \label{eq:CGL_Galerkin}
dX=(i+\gamma) \left(HX-S_N( |S_N X|^2 S_N X)\right)dt +\sqrt{2\gamma}dW. 
\end{equation}
It follows from Proposition \ref{prop:prop4} that  
$X^N$ is globally defined on $L^p(\R)$ if $p\ge 3$,  
and depends  continuously on the initial data $X_0$. 
We remark that the results in Section 2 also lead 
to the convergence of $X^N$ to $X$ in $C([0,T], L^p)$ for any $T<T^*$. 
Consider now 
the transition semigroup for the equation (\ref{eq:CGL_Galerkin}) and denote it by $P_t^N$ for $t>0$. 
Notice that $\rho_N$ is an invariant measure for the transition semigroup $P_t^N$ for the same reason as in Proposition \ref{invariance}. 
\vspace{3mm}

\noindent
{\bf Proof of Theorem \ref{thm:thm1}}. 
Let any $T>0$. We wish to prove that there exists a constant $C_T$ such that 
\begin{equation} \label{eq:global_X}
 \int_{L^p} \E\Big(\sup_{t\in [0,T^* \wedge T]} |X(t, y)|_{L^p}\Big) d\rho(y) \le C_T. 
\end{equation}
If this bound is true, we may say that there exists a $\rho$-measurable set $\mathcal{O} \subset L^p$ such that 
$\rho(\mathcal{O})=1$, and for each $y \in \mathcal{O}$, $\sup_{t\in [0,T^* \wedge T)} |X(t, y)|_{L^p}<+\infty$ a.s.;
the fixed point argument then implies that $T^*(y)\ge T$, a.s. 
\vspace{3mm}

In order to derive rigorously the estimate (\ref{eq:global_X}), 
we will use the approximation $X^N$. Thus the first step here is to show that  
there exists a constant $C_T$ which does not depend on $N$, such that
\begin{eqnarray} \label{eq:global_XN}
\int_{L^p} \E\Big(\sup_{t\in [0,T]} |X^N(t, y)|_{L^p}\Big) d\rho_N (y) \le C_T.  
\end{eqnarray}
We write the equation of $X^N$ in the mild form. 
\begin{eqnarray*}
X^N(t, X_0)&=&e^{t(i+\gamma)H}X_0-e^{t(i+\gamma)H}Z_{\infty}(0) \\
&& -(i+\gamma)\int_0^t e^{(t-s)(i+\gamma)H} S_N(|S_NX^N|^2 S_NX^N)(s) ds +Z_{\infty}(t).
\end{eqnarray*}
Taking the supremum in time on $[0,T]$, using Lemma \ref{lem:linear_deterministic} 
and the boundness of $S_N$ in $L^p(\R) \cap L^{3p}(\R)$, 
we have 
\begin{eqnarray*}
\sup_{t\in [0,T]} |X^N(t, X_0)|_{L^p}  
&\le& |X_0|_{L^p} +  2 \sup_{t\in [0,T]} |Z_{\infty}(t)|_{L^p} 
+ C_{\gamma}\int_0^T | S_N(|S_N X^N(s)|^2 S_N X^N(s))|_{L^p} ds\\
&\le& |X_0|_{L^p} + 2 \sup_{t\in [0,T]} |Z_{\infty}(t)|_{L^p} 
+ C_{\gamma} \int_0^T |S_N X^N(s)|_{L^{3p}}^3 ds\\
 &\le& |X_0|_{L^p} + 2 \sup_{t\in [0,T]} |Z_{\infty}(t)|_{L^p} 
+  C_{\gamma} \int_0^T |X^N(s)|_{L^{3p}}^3 ds.
\end{eqnarray*}
We then take the expectation, and we have 
\begin{eqnarray*}
\E\Big(\sup_{t\in [0,T]} |X^N(t, X_0)|_{L^p} \Big) 
&\le & |X_0|_{L^p} + 2 \E\left(|Z_{\infty}|_{L^{\infty}(0,T;L^p)}\right) + C_{\gamma} \E \int_0^T |X^N (s)|_{L^{3p}}^3 ds.\\
\end{eqnarray*}
The second term is bounded by a constant $M_{p,T}$ which does not depend on $X_0$ thanks to Lemma \ref{regularityZ}. 
Integrating in $\rho_N$ thus gives, 
\begin{eqnarray*}
&& \int_{L^p} \E\Big(\sup_{t\in [0,T]} |X^N(t, X_0)|_{L^p} \Big) d\rho_N (X_0)\\
&\le & \int_{L^p} |X_0|_{L^p} d\rho_N (X_0) + M_{p,T}
+ C_{\gamma} \int_{L^p} \int_0^T \E |X^N(s)|_{L^{3p}}^3 ds \,d\rho_N(X_0). 
\end{eqnarray*}
The third term in the right hand side is written as follows :
\begin{equation*}
\int_0^T \E \int_{L^p} |X^N(s, X_0)|_{L^{3p}}^3 ds \,d\rho_N(X_0)=T\int_{L^p} |X_0|_{L^{3p}}^3 d\rho_N(X_0);  
\end{equation*}
indeed, let $\Phi^m(v)=|v|_{L^{3p}}^3 \wedge m$ for any $v \in L^p$, which is a bounded Borel function on $L^p$. 
\begin{eqnarray*}
\int_0^T \E \int_{L^p} |X^N(s,X_0)|_{L^{3p}}^3 ds\, d\rho_N(X_0) 
&=& \int_0^T \lim_{m\to \infty} \E \int_{L^p} \Phi^m (X^N(s,X_0)) d\rho_N(X_0) \, ds\\
&=& \int_0^T \lim_{m\to \infty}  \int_{L^p} P_s^N \Phi^m(X_0) d\rho_N(X_0) \, ds \\
&=& \int_0^T \lim_{m\to \infty} \int_{L^p} \Phi^m(X_0) d\rho_N (X_0) \, ds\\
&=& T \int_{L^p} |X_0|_{L^{3p}}^3 d\rho_N (X_0).  
\end{eqnarray*}
Here, we have used the invariance of $\rho_N$ in the third equality. 
We shall show that 
\begin{equation*}
\int_{L^p} |X_0|_{L^{3p}}^3 d\rho_N(X_0) 
\end{equation*}
is bounded independently of $N$. This can be seen as follows.  
Recalling the definition \eqref{defrhoN}, noting that for all $t\ge 0$, $\mathcal{L}(Z_{\infty}(t))=\mathcal{L}(Z_{\infty}(0))=\mu$,  
with $\mathrm{supp} \mu \subset \bigcap_{p>2} L^p$, and 
Lemma \ref{regularityZ}, we have 
\begin{eqnarray*}
\int_{L^p} |X_0|_{L^{3p}}^3 d\rho_N(X_0)
&=& \Gamma_N^{-1} \int_{L^p} |X_0|_{L^{3p}}^3 e^{-\frac{1}{4}|S_N X_0|_{L^4}^4} d\mu(X_0) \\
&\le&\Gamma_N^{-1} \, \E(|Z_{\infty}(0)|_{L^{3p}}^3) \le \Gamma_N^{-1} M_{3p,T}^3.
\end{eqnarray*} 
As for the normalizing constant $\Gamma_N,$ 
the bound $|S_N X_0|_{L^4} \le C |X_0|_{L^4}$ with a constant $C$ independent of $N$ implies 
\begin{equation} \label{bound:normalizedcst}
\Gamma_N= \int_{L^p} e^{-\frac{1}{4}|S_N X_0|_{L^4}^4} d\mu(X_0) \ge \int_{L^p} e^{-\frac{C}{4}|X_0|_{L^4}^4} d\mu(X_0) >0,  
\end{equation}
i.e. $\Gamma_N$ is bounded from below independently of $N.$ 

For the same reason, the first term is bounded independently of $N$ since 
\begin{eqnarray*}
\int_{L^p} |X_0|_{L^p} d\rho_N(X_0) \le \Gamma_N^{-1} M_{p,T},  
\end{eqnarray*}
with (\ref{bound:normalizedcst}). This concludes the bound (\ref{eq:global_XN}).
\vspace{3mm}

Recall that if $p\ge 3$, then by Proposition \ref{prop:prop4}, 
$X^N(\cdot,X_0) \to X(\cdot, X_0)$ in $C([0,T\wedge T^*), L^p(\R))$ as $N \to \infty$ for any $X_0 \in L^p \cap L^4$.  
Since 
$e^{-\frac{1}{4}|S_N X_0|_{L^4}^4} \to e^{-\frac{1}{4}|X_0|_{L^4}^4}$ as $N \to \infty$ for any $X_0 \in L^p \cap L^4$, 
we obtain by Fatou's lemma
\begin{eqnarray*}
\int_{L^p} \E\Big(\sup_{t\in [0,T \wedge T^*]} |X(t, X_0)|_{L^p}\Big) d\rho(X_0) 
&=&
\int_{L^p} \E\Big(\sup_{t\in [0,T \wedge T^*]} |X(t, X_0)|_{L^p}\Big) e^{-\frac{1}{4}|X_0|_{L^4}^4} d\mu(X_0) \\
&\le& 
\liminf_{N\to \infty} 
 \int_{L^p} \E\Big(\sup_{t\in [0,T\wedge T^*]} |X^N(t, X_0)|_{L^p} \Big) e^{-\frac{1}{4}|S_N X_0|_{L^4}^4} d\mu(X_0) \\
&=&
 \liminf_{N\to \infty} 
 \int_{L^p} \E\Big(\sup_{t\in [0,T\wedge T^*]} |X^N(t, X_0)|_{L^p} \Big) d\rho_N (X_0) \\
&\le& C_T, 
\end{eqnarray*}
which concludes (\ref{eq:global_X}). 

In conclusion, for a fixed $T>0$, there exists a $\rho$ - measurable set $\mathcal{O}_T \subset L^p(\R)$ for any $p\ge 3$ such that 
$\rho(\mathcal{O}_T)=1$ and that for $X_0 \in \mathcal{O}_T,$ the solution $X$ exists a.s. up to time $T$.
For each $T_n$, $n\in \N$ such that $T_n \to +\infty$ as $n\to \infty$, consider 
$\mathcal{O}_{T_n}$ and set $\bar{\mathcal{O}}:= \bigcap_{n} \mathcal{O}_{T_n}$. 
Then, $\rho(\bar{\mathcal{O}})=1$. The proof is completed.  
\hfill\qed
\vspace{3mm}

Now, we can define the transition semi-group $(P_t)_{t\ge 0}$ associated with \eqref{eq:CGL} :
for all $t\ge 0$ and any $\varphi \in C_b(L^p, \R)$ with $p\ge 3$, and for $\rho$ - a.e. $y$ (i.e., for $y\in \bar{\mathcal{O}}$), let
$$P_t \varphi (y)=\E (\varphi(X(t, y))).$$
We end this section with the proof of the invariance of the measure $\rho$ for the semi-group $(P_t)_{t\ge 0}$.
Let $\varphi \in C_b(L^p, \R)$ and let us prove that
\begin{equation}
\label{invariancerho}
\int_{L^p} \varphi (y) d\rho(y) =\int_{L^p} P_t \varphi(y) d\rho (y). 
\end{equation} 
First, by the invariance of $\rho_N$, we have 
\begin{equation*}
\int_{L^p} \varphi (y) d\rho_N (y) =\int_{L^p} P_t^N \varphi(y) d\rho_N (y). 
\end{equation*}
Next,
\begin{eqnarray*}
\Big|\int_{L^p} \varphi (y) d\rho_N (y)- \int_{L^p} \varphi (y) d\rho(y) \Big| 
&\le& \Big|\Gamma_N^{-1} \Big(\int_{L^p} \varphi (y) e^{-\frac{1}{4}|S_N y|_{L^4}^4}d\mu (y)
- \int_{L^p} \varphi (y) e^{-\frac{1}{4}|y|_{L^4}^4} d\mu(y)\Big) \Big| \\
&& + \Big|(\Gamma_N^{-1}-\Gamma^{-1}) \int_{L^p} \varphi (y) e^{-\frac{1}{4}|y|_{L^4}^4} d \mu (y) \Big| \\
&\le& C(1+|\phi|_{L^{\infty}})|\Gamma-\Gamma_N|,\\
\end{eqnarray*}
where we have used \eqref{bound:normalizedcst}, and the right hand side above tends to zero by the dominated convergence
theorem. On the other hand,
\begin{eqnarray*}
\Big|\int_{L^p} P^N_t \varphi (y) d\rho_N (y)- \int_{L^p} P_t \varphi (y) d\rho(y) \Big| 
&\le& \Big|\int_{L^p} \E (\varphi(X(t,y))) d\rho_N(y)- \int_{L^p} \E (\varphi(X(t,y)))d\rho(y)\Big| \\
&& + \Big|\int_{L^p} [\E(\varphi(X^N(t,y)))-\E (\varphi(X(t, y)))] d\rho_N (y)\Big|\\
&\le& |\varphi|_{L^{\infty}} \int_{L^p} \big| e^{-\frac14 |S_Ny|^4_{L^4}}-e^{-\frac14|y|^4_{L^4}}\big| d\mu(y) \\
&& + \int_{L^p} |\E(\varphi(X^N(t,y)))-\E (\varphi(X(t, y)))| d\mu(y),
\end{eqnarray*}
and we conclude thanks to the dominated convergence theorem and the fact 
that $X^N(\cdot, y)$ converges to $X(\cdot, y)$ a.s. in $C([0,T], L^p)$ for any $y \in \bar{\mathcal{O}}$. 
\hfill\qed

\section{Strong Feller property and global existence for all initial data}

We first prove the strong Feller property for the semigroups $(P_t^N)_{t\ge 0}$, associated with equation 
\eqref{eq:CGL_Galerkin}, uniformly in $N$. We denote by 
$\mathcal B_b(E)$ the space of borelian bounded real valued functions on a Banach space $E$ and, for
$\varphi\in \mathcal B_b(E)$, $\|\varphi\|_0=\sup_{x\in E}|\varphi(x)|$.
\begin{prop}\label{p8} 
Let $p\ge 3$ and $N\in\N$. For any $\varphi\in \mathcal B_b(L^p(\R))$, and any $t>0$,
$P^N_t\varphi$ is a continuous function on $L^p(\R)$. Moreover,  for any $T>0$, there exists a constant $c_\gamma(T)$, independent of 
$N$, such that for any $X_0 \in \bar{\mathcal{O}}$, $h\in L^p(\R)$ such that $X_0 +h \in \bar{\mathcal{O}}$, and any $T>0$,
\begin{equation}\label{e7.1}
|P_T^N\varphi(X_0+h)-P_T^N\varphi(X_0)|\le  c_\gamma(T)\|\varphi\|_0 \left(T^{-1+\frac1p}+1\right)\left(|X_0|_{L^p(\R)}+1\right)^2|h|_{L^p(\R)}.
\end{equation}
\end{prop}

\proof  We use a classical coupling argument, based on a control problem and Girsanov Theorem
 (see e.g. \cite{priola06}). Let $p\ge 3$, $T>0$, $h\in L^p(\R)$ and define for $t\in [0,T]$ :
$$
k(t)= - \frac{1}{T} e^{(i+\gamma) H t} h,
$$
and
$$
\bar h(t)= \frac{T-t}{T} e^{(i+\gamma)H t}h.
$$
Then 
$$
\frac{d}{dt}\bar h=(i+\gamma)H \bar h + k,\quad \bar h (0)=h,\quad \bar h(T)=0.
$$
Let $\theta\in C_0^\infty(\R)$ be a cut-off function such that $\theta(t)=1$ for $t\in [0,1]$ and  $\theta(t)=0$ for $t\ge 2$. Define,
for $y\in L^p(\R)$, 
$F_R(y)=\theta\left(\frac{|S_Ny|_{L^p(\R)}^p}R\right)S_N( |S_N y|^2 S_N y)$.
We consider the truncated version of \eqref{eq:CGL_Galerkin} :
\begin{equation} \label{eq:galcut}
dX= (i+\gamma) \left(HX- F_R(X)\right)dt +\sqrt{2\gamma}dW
\end{equation}
with initial data $X(0)=X_0$. We denote by $X_R(\cdot,X_0)$ its solution, whose existence and uniqueness in 
$C([0,T];L^p(\R))$ follows from the arguments of Section 3, and the Lipschitz property of $F_R$ (note that 
$F_R(y)=0$ if $|S_N y|_{L^p}^p\ge 2R$).
For $\eps >0$, define $Y^{\eps}(t)= X_R(t,X_0)+\eps \bar h(t)$. Then
$$
dY^{\eps}= (i+\gamma) \left(H Y^{\eps}-F_R(Y^{\eps})\right)dt + \sqrt{2\gamma}\, (G_\eps dt +dW)
$$
with 
$$
G_\eps(t)= \frac{(i+\gamma)}{\sqrt{2\gamma}} \left[F_R(Y^{\eps}(t))-F_R(X_R(t,X_0))\right] + \eps k(t).
$$
Since $F_R$ is bounded and $k\in L^2(0,T;L^p(\R))$, we may use Girsanov's transform and deduce  for each $\varphi\in C_b^1(L^p(\R))$
$$
\E[\varphi(Y^\eps(T))]=\E\left[\varphi(X_R(T,y))\rho_\eps(T)\right],
$$
where
$$
\rho_\eps(T)=\exp\left\{-\gamma\int_0^T|G_\eps(t)|_{L^2(\R)}^2\,dt - \sqrt{2\gamma}\int_0^T( G_\eps(t),dW(t))_{L^2(\R)}\right\},
$$
and $y=X_0+\eps h$.
Clearly $Y^\eps(T)=X_R(T,X_0)$ and
we obtain 
\begin{equation}\label{e7.2}
\E[\varphi(X_R(T,X_0))]= \E[\varphi(X_R(T,X_0+\eps h))\rho_\eps(T)].
\end{equation}
On the other hand, noting that for $\eps=0$ we have $G_\eps=0$ and $\rho_\eps=1$,
$$
D_\eps G_\eps(t)|_{\eps=0}=k(t)+ \frac{(i+\gamma)}{\sqrt{2\gamma}} F'_R(X_R(t,X_0))\cdot \bar h(t)
$$
and
$$
D_\eps \rho_\eps(T)|_{\eps=0}
=-\sqrt{2\gamma} \int_0^T\Big( k(t)+\frac{(i+\gamma)}{\sqrt{2\gamma}} F_R'(X_R(t,X_0))\cdot \bar h(t),dW(t)\Big)_{L^2(\R)}.
$$
Since $\varphi\in C^1_b(L^p(\R))$ , $\bar h, k\in L^2(0,T;L^p(\R))$, we  may differentiate \eqref{e7.2} with respect to $\eps$ and take $\eps=0$. we obtain :
\begin{eqnarray*}
& & \E[D\varphi(X_R(T,X_0))\cdot \left(DX_R(T,X_0)\cdot h\right) ]\\
& = &\sqrt{2\gamma} \,  \E\left[\varphi(X_R(T,X_0))\int_0^T\Big( k(t)+ \frac{(i+\gamma)}{\sqrt{2\gamma}}
F_R'(X_R(t,X_0))\cdot \bar h(t),dW(t)\Big)_{L^2(\R)}\right].
\end{eqnarray*}
By the chain rule, the left hand side above is equal to 
$$
D\left[ \E\left(\varphi(X_R(T,X_0))\right)\right]\cdot h= DP_T^R\varphi(X_0)\cdot h,
$$
 where $(P^R_t)_{t\ge 0}$ is the transition semigroup associated to equation \eqref{eq:galcut}. Note that $P^R_t\varphi(X_0)$
 tends to $P^N_t\varphi(X_0)$ when  $R$ goes to infinity, for any $X_0 \in L^p(\R)$.
 On the other hand, the right hand side may be bounded using the It\^o isometry and the Cauchy Schwarz inequality, to get
\begin{equation} \label{estimrhs}
DP_T^R\varphi(X_0)\cdot h \le  C_\gamma \|\varphi\|_0 \left(\int_0^T |k(t)|_{L^2(\R)}^2+ \E\, |F_R'(X_R(t,X_0))\cdot \bar h(t)|_{L^2(\R)}^2dt
\right)^{1/2}
\end{equation}
We need the following result.
\begin{lem}\label{l7.1}
Let $p\in [2,\infty]$, then, for $t>0$, $e^{(i+\gamma)tH}$ maps $L^p(\R)$ into $L^2(\R)$ and there exists a constant $C_{\gamma,p}$ such that
for all $t>0$, and all $f\in L^p(\R)$,
$$
|e^{(i+\gamma)tH}f|_{L^2(\R)}\le C_{\gamma,p} \, t^{-\frac12+\frac1{p}} |f|_{L^p(\R)}.
$$
\end{lem}
\proof
Let $f\in L^\infty(\R)$ and $u=(-H)^{-1}f$, then 
$
-\partial^2_x u+ x^2 u= f
$
and 
\begin{eqnarray*}
 |(-H)^{-1/2}f|_{L^2(\R)}^2&=&|(-H)^{1/2}u|_{L^2(\R)}^2=\int_\R |\partial_x u(x)|^2 + x^2 u^2(x) dx\\
&=& \int_\R f(x)u(x)dx  \le |f|_{L^\infty}\int_\R |u(x)|dx.
\end{eqnarray*}
Note that 
$$
\int_\R |u(x)|dx\le \left(\int_\R \frac1{1+x^2}dx\right)^{1/2} \left(|u|_{L^2(\R)}^2+|xu|_{L^2(\R)}^2\right)^{1/2}
= \sqrt{\pi} |(-H)^{1/2}u|_{L^2(\R)},
$$
and it follows :
\begin{equation} \label{estiminf}
|(-H)^{-1/2}f|_{L^2(\R)}\le \sqrt{\pi}|f|_{L^\infty}.
\end{equation}
On the other hand, writing $v=\sum_{k=0}^{\infty} (v,e_k)_{L^2} e_k$, we have for any $v\in L^2(\R)$,
$$
|(-H)^{1/2}e^{(i+\gamma)tH} v|_{L^2(\R)}^2 = \sum_{k=0}^{\infty} \lambda_k^2 |e^{-(i+\gamma)t\lambda_k^2}|^2(v,e_k)_{L^2}^2\\
\le \frac{C_\gamma}{t} |v|_{L^2}^2,
$$
so that \eqref{estiminf} implies
\begin{eqnarray*}
|e^{(i+\gamma)tH}f|_{L^2(\R)}&= &|(-H)^{1/2}e^{(i+\gamma) tH}(-H)^{-1/2}f|_{L^2(\R)}\le C_\gamma t^{-1/2} |(-H)^{-1/2}f|_{L^2(\R)}\\
&\le & \sqrt{\pi} C_\gamma t^{-1/2} |f|_{L^{\infty}(\R)}.
\end{eqnarray*}
This is the result for $p=\infty$. For $p=2$, the result is clear and the general case follows by interpolation.
\hfill$\qed$
\medskip

\begin{rem}
This result is not optimal. The exponent in $t$ can actually be taken as $-\beta$ for any $\beta>\frac14-\frac1{2p}$.
\end{rem}
To end the proof of Proposition \ref{p8}, we need to bound the right hand side of \eqref{estimrhs}. First, we write
\begin{eqnarray} 
\nonumber \int_0^T |k(t)|_{L^2(\R)}^2dt &= & \frac{1}{T^2} \int_0^T|e^{(i+\gamma) Ht} h|^2_{L^2(\R)} dt\\
\nonumber & \le &\frac{c}{T^2}\int_0^T t^{-1+\frac2p}|h|_{L^p(\R)}^2dt\\  \label{estimk}
& \le &c T^{\frac2p-2}|h|_{L^p(\R)}^2, 
\end{eqnarray}
thanks to Lemma \ref{l7.1}. Next, we compute, for $x, y \in L^p(\R)$,
\begin{eqnarray*}
F'_R(x)\cdot y&=& \frac{p}{R}\theta' \left( \frac{|S_N(x)|_{L^p}^p}{R}\right) \mathcal{R}e \left(|S_Nx|^{p-2}S_Nx,S_Ny\right)_{L^2}
S_N\left(|S_Nx|^2 S_Nx\right)\\
& & + \, \theta\left(\frac{|S_N(x)|_{L^p}^p}{R}\right) S_N\left(|S_Nx|^2S_Ny\right) \\
& & +\, \theta\left(\frac{|S_N(x)|_{L^p}^p}{R}\right) S_N\left(2 \mathcal{R}e(S_Nx\cdot\overline{S_Ny})S_Nx\right).
\end{eqnarray*}
Hence, the second term in \eqref{estimrhs} is bounded as follows, using H\"older inequality and Lemma \ref{lem:linear_deterministic} :
\begin{eqnarray}
\nonumber \E\int_0^T |F_R'(X_R(t,X_0))\cdot \bar h(t)|_{L^2(\R)}^2dt & \le
& \frac{C}{R^2} 
\E \int_0^T |S_N(X_R(t,X_0)|^{2(p-1)}_{L^p} |\bar h(t) |_{L^p}^2 |S_N(X_R(t,X_0))|_{L^6}^6 dt\\
\nonumber & & + \, C \, \E\int_0^T |S_N(X_R(t,X_0))|_{L^4}^4 |\bar h(t)|_{L^{\infty}}^2 dt \\
\nonumber & \le & \frac{C}{R^{2/p}} \E \int_0^T |S_N(X_R(t,X_0))|_{L^6}^6 |h|_{L^p}^2 dt\\
& & +\, C \, \E \int_0^T t^{-\frac{1}{p}} |S_N(X_R(t,X_0))|_{L^4}^4 |h|_{L^p}^2 dt. \label{estimF1}
\end{eqnarray}
Next, we decompose $X_R(t,X_0)=v_N^R(t,0) +Z(t)$, with
\begin{equation} \label{defZ}
Z(t) =e^{(i+\gamma)tH} X_0+\sqrt{2\gamma} \int_0^t e^{(i+\gamma)(t-s)H}dW(s).
\end{equation}
If $3\le p\le 4$, we use Lemma \ref{lem:linear_deterministic} to get $|e^{(i+\gamma)tH}X_0|_{L^4}^4\le C_\gamma t^{-2/p+1/2} |X_0|_{L^p}^4$;
if $p>4$, we interpolate the inequality of Lemma \eqref{l7.1} (with $p=\infty$) and the boundedness of the operator $e^{(i+\gamma)tH}$
in $L^q(\R)$ for any $q$ with $4<q<p$ to get $|e^{(i+\gamma)tH}X_0|_{L^4}^4\le C_\gamma t^{4/p-1} |X_0|_{L^p}^4$.
In both case, this and Lemma \ref{regularityZ} imply, since $Z$ is Gaussian, and for any integer $m$,
\begin{equation}
\label{estimZ}
\E |Z|_{L^4(0,T;L^4(\R))}^{4m} \le c \left[\int_0^T \big|e^{(i+\gamma)tH}X_0\big|_{L^4(\R)}^4 dt \right]^m+M_{\gamma,p,T,m} \le C_{\gamma,p,T} (|X_0|_{L^p(\R)}+1)^{4m}.
\end{equation}
Since moreover the energy inequality \eqref{ineq:energy} (with $\eta=0$), which is easily seen to hold also for
the cut-off equation satisfied by $v_N^R$, implies
$$
\sup_{t\in[0,T]}|v_N^R(t)|_{L^2(\R)}^2 \le C_{\gamma,T} |Z|_{L^4(0,T;L^4(\R))}^4
$$
and
$$
\int_0^T |\partial_x v_N^R(t)|_{L^2(\R)}^2 dt \le C_{\gamma,T} |Z|_{L^4(0,T;L^4(\R))}^4,
$$
we deduce from \eqref{estimZ} and the Gagliardo-Nirenberg inequalities :
$$
|v_N^R(t)|_{L^4(\R)}^4\le C |v_N^R(t)|_{L^2(\R)}^3 |\partial_x v_N^R(t)|_{L^2(\R)}
$$
and
$$
|v_N^R(t)|_{L^6(\R)}^6\le C |v_N^R(t)|_{L^2(\R)}^4 |\partial_x v_N^R(t)|_{L^2(\R)}^2,
$$
that
$$
\E |v_N^R|_{L^8(0,T;L^4(\R))}^4 \le C_{\gamma,T} \E |Z|_{L^4(0,T;L^4(\R))}^8 \le C_{\gamma,p,T} (|X_0|_{L^p(\R)}+1)^8,
$$
and 
$$
\E |v_N^R|_{L^6(0,T;L^6(\R))}^6 \le C_{\gamma,T} \E |Z|_{L^4(0,T;L^4(\R))}^{12} \le C_{\gamma,p,T} (|X_0|_{L^p(\R)}+1)^{12}.
$$
Plugging these inequalities into the right hand side of \eqref{estimF1}, after using Cauchy-Schwarz inequality, shows that
\begin{equation} \label{estimF2}
\E\int_0^T |F_R'(X_R(t,X_0))\cdot \bar h(t)|_{L^2(\R)}^2dt \le C_{\gamma,p,T} (|X_0|_{L^p}+1)^8\big[1+\frac{1}{R^{2/p}}(|X_0|+1)^4\big]|h|_{L^p(\R)}^2.
\end{equation}
Finally, \eqref{estimrhs}, \eqref{estimk}, and \eqref{estimF2} imply after taking the limit $R\to \infty$ :
$$
DP_T^N\varphi(X_0)\cdot h \le  C_{\gamma,T,p}\|\varphi\|_0 \left(|X_0|_{L^p(\R)}+1\right)^4|h|_{L^p(\R)}.
$$
It remains  only  to approximate $\varphi\in \mathcal B_b(L^p(\R))$ by a sequence $(\varphi_n)_{n\in\N}$ in $C^1_b(L^p(\R))$,
which converges pointwise to $\varphi$ and are uniformly bounded,  to end the proof.
\hfill$\qed$
\medskip

Now, in order to prove Theorem \ref{t3}, we use the fact that there exists a set $\mathcal O\subset L^p(\R)$ such that $\rho( \mathcal O)=1$ and for any 
$X_0\in \mathcal O$ we have $\sup_{t\in[0,T]} |X(t,X_0)|_{L^p(\R)} <\infty$ a.s. Note that this holds in any probability
space since we have constructed strong solutions.
Letting $N\to\infty$ in \eqref{e7.1} for $X_0,\; X_0+h \in \mathcal O$, we deduce that $(P_t)_{t\ge 0}$ can be extended 
uniquely into a semigroup on $L^p(\R;\C)$ which has the strong Feller property. This is possible since we know that the 
support of $\rho$ in the topology of $L^p(\R)$ is $L^p(\R)$.

Let us now choose $X_0\in L^p(\R)$ and take a sequence $(X_{0,n})_{n\in\N}$ in ${ \mathcal O}$ such that 
$X_{0,n}\to X_0$ in $L^p(\R)$. 
Considering the corresponding solution $X^{N}(t,X_{0,n})$ of \eqref{eq:CGL_Galerkin},
we use again the splitting: $X^{N}(t,X_{0,n})=v_{N,n}+Z_n$, with $v_{N,n}$ satisfying \eqref{eq:deterministic_Galerkin1} with 
$$z(t)=Z_n(t):=e^{(i+\gamma)tH}X_{0,n}+\sqrt{2\gamma}\int_0^t e^{(t-s)(i+\gamma)H}dW(s).$$
Clearly, 
$Z_n\to Z$ in $C([0,T]; L^p(\R))$ a.s. where $Z$ is defined in \eqref{defZ}.
Moreover, using  \eqref{ineq:energy} for  $v_{N,n}$ and letting
$N\to \infty$, we deduce that $v_n=X(\cdot,X_{0,n})-Z_n$ is bounded uniformly in $n$ in 
$L^2(\Omega; L^\infty(0,T,L^2(\R))\cap L^2(0,T;{\mathcal W}^{1,2}(\R) ))$$\cap L^4(\Omega; L^4(0,T;L^4(\R)))$. 
By standard arguments, using again equation \eqref{eq:deterministic_Galerkin1},
it follows that, for any negative $s$, the 
sequence of laws of 
$(X(\cdot,X_{0,n}))_{n\in\N}$ is tight in $C([0,T];{\mathcal W}^{s,2}(\R))\cap L^2(0,T;{\mathcal W}^{1+s,2}(\R))
\cap C^1([0,T];{\mathcal W}^{s-2,2}(\R) )+C([0,T]; L^p(\R))$. It is then standard to check that any limit point is the law of a 
martingale solution to \eqref{eq:CGL}. Let us denote by $\tilde X(\cdot,X_0)$ such a solution, which is defined
on a probability space $(\tilde \Omega, \tilde{ \mathcal F}, \tilde{\mathbb P})$.
Take a function $\varphi\in C_b(L^p(\R))\cap C_b({\mathcal W}^{s,2}(\R))$ ; then we can write for any $t>0$
$$P_t\varphi(X_0)=\lim_{n\to \infty} P_t\varphi(X_{0,n})=\lim_{n\to\infty}\E(\varphi(X(t,X_{0,n})))=
\tilde\E(\varphi(\tilde X(t,X_{0})).$$
This proves that the law of $\tilde X(t,X_{0})$ is in fact $P_t^*\delta_{X_0}$. 

Now, we make use of the following result, which is proved as in \cite{DPD03} Lemma 7.7 (see also \cite{flandoli97} 
for the original idea).
\begin{prop}\label{p9}
The semigroup $(P_t)_{t\ge 0}$ is irreducible on $L^p(\R)$.
\end{prop}

This and Proposition \ref{p8} mean that $(P_t)_{t> 0}$ is regular and that all the measures $P^*_t\delta_{X_0}$ are 
equivalent and are also equivalent to the invariant measure $\rho$, which is the unique invariant measure of $(P_t)_{t\ge 0}$
(see \cite{dpz0}, chapter 4). It follows that for any $X_0\in L^p(\R)$ and any $t_0>0$,
$$
{\tilde {\mathbb P}}(\tilde X(t_0,X_0)\in \mathcal O)=\rho(\mathcal O)=1
$$
so that
$$
\sup_{t \in [t_0,T]}|\tilde X(t,X_0)|_{L^p(\R)}<\infty, \quad {\tilde {\mathbb P}} \; a.s.
$$
By the construction of local solutions, we know that there exists a stopping time $\tau^*(X_0)$ with
$0< \tau^*(X_0) \le T$, a.s. such that 
$$
\sup_{t \in [0,\tau^*(X_0))}|\tilde X(t,X_0)|_{L^p(\R)}<\infty, \quad {\tilde {\mathbb P}}\; a.s.
$$
It follows, since $\tau^*(X_0)>0$, a.s., that
$$
\sup_{t\in [0,T]}|\tilde X(t,X_0)|_{L^p(\R)}<\infty,
\quad {\tilde {\mathbb P}}\; a.s.,$$
and we have constructed a global solution in $L^p(\R)$ for any $X_0\in L^p(\R)$. This solution is clearly pathwise unique
in $C([0,T];L^p(\R))$ and either using Yamada-Watanabe theorem or Gyongy, Krylov method we may deduce 
global existence and uniqueness of a strong solution (in the probabilistic sense) in any probabilistic space.
\hfill $\qed$
\medskip

\section{Convergence to equilibrium}


We recall that $X_N(.,y_N)$ is the solution of (\ref{eq:CGL_finite_Galerkin}) with initial data $X_N(0)=:y_N\in E_N$. 
By  Proposition \ref{invariance}, the measure $\tilde{\rho}_N$ is invariant by the flow of $X_N$, and by Proposition \ref{prop:poincare},
it satisfies the Poincar\'e inequality (\ref{ineq:poincare}).
Let $(\tilde{P}_t^N)_{t\ge 0}$ be the transition semigroup corresponding to (\ref{eq:CGL_finite_Galerkin}). 
Using the above properties, we have

\begin{prop} \label{prop:prop7}
Let $\phi_N \in C_b (E_N, \R)$. 
Then,  for $t \ge 0$ and $\gamma >0$, 
\begin{equation*}
\frac{d}{dt}\int_{E_N} |\tilde{P}_t^N \phi_N(y)-\bar{\phi}_N^{\tilde{\rho}_N}|^2 d\tilde{\rho}_N(y) 
+2\gamma \int_{E_N} |\tilde{P}_t^N \phi_N(y)-\bar{\phi}_N^{\tilde{\rho}_N}|^2 d\tilde{\rho}_N(y)
\le  0, 
\end{equation*}
i.e., 
$$\int_{E_N} |\tilde{P}_t^N \phi_N(y)-\bar{\phi}_N^{\tilde{\rho}_N}|^2 d\tilde{\rho}_N(y) \le e^{-2\gamma t}
\int_{E_N} |\phi_N(y)-\bar{\phi}_N^{\tilde{\rho}_N}|^2 d\tilde{\rho}_N(y),$$
where $\bar{\phi}_N^{\tilde{\rho}_N}=\int_{E_N} \phi_N(y) d\tilde{\rho}_N (y).$
\end{prop}

\proof Using a density argument, we may assume $\bar{\phi}_N^{\tilde{\rho}_N}=0$, and $\phi_N \in C_b^1 (E_N, \R)$.
Let $\mathcal{L}_N$ be the generator of $\tilde{P}_t^N$, associated with \eqref{eq:CGL_finite_Galerkin}. 
As was seen in the proof of Proposition \ref{prop:poincare} for
the generator $\mathcal{M}_N$, 
using the Kolmogorov equation and the invariance of 
the measure $\tilde{\rho}_N$, 
\begin{equation*} 
0=\int_{E_N} \mathcal{L}_N (\tilde{P}_t^N \phi_N)^2  d\tilde{\rho}_N 
=\frac{d}{dt} \int_{E_N}|\tilde{P}_t^N \phi_N|^2 d\tilde{\rho}_N + 2\gamma\int_{E_N}|D\tilde{P}_t^N \phi_N|_{E_N}^2 
d\tilde{\rho}_N. 
\end{equation*}
The use of inequality (\ref{ineq:poincare}) implies the desired result. \hfill\qed
\medskip

Now taking the limit $N\to \infty$, we shall show the following proposition. 

\begin{prop} \label{prop:prop6}
Let $p\ge 3$ and let $\phi \in C_b(L^p(\R);\R).$ Then there exists $\phi_N \in C_b(E_N, \R)$ such that  for all $t\ge 0$,
\begin{equation*}
\lim_{N \to \infty} \int_{E_N} |\tilde{P}_t^N \phi_N(y)-\bar{\phi}_N^{\tilde{\rho}_N}|^2  d\tilde{\rho}_N (y) 
= \int_{L^p} |P_t \phi(y) -\bar{\phi}^{\rho}|^2  d\rho(y),~ \mbox{and}
~ \lim_{N\to \infty} \bar{\phi}_N^{\tilde{\rho}_N}=\bar{\phi}^{\rho},
\end{equation*}
where $\bar{\phi}^{\rho}=\int_{L^p} \phi(y) d\rho (y).$
\end{prop}

Note that the proof of Theorem \ref{thm:thm2} follows immediately from Proposition \ref{prop:prop7} and Proposition \ref{prop:prop6}.
The proof of Proposition \ref{prop:prop6} makes use of the next Lemma.

\begin{lem} \label{lem:cylindrical_approx}
Let $1 \le p \le +\infty.$ 
Let $\psi \in C_b (L^p(\R); \R)$ and let $\psi _{|_{E_N}} (u):=\psi(S_N(u))$ for $u\in L^p(\R)$.  
Then, for any $q$ with $1 \le q < +\infty$, $\big|\psi_{|_{E_N}} -\psi \big|_{L^q(L^p, d\rho)}$ tends to $0$ as $N$
goes to infinity.
\end{lem}

\proof Note that $S_N$ is continuous from $L^p(\R)$ to $L^p(\R)  $, so
$\psi_{|_{E_N}}$ is continuous from $L^p(\R)$ to $\R$.
Also $\psi_{|_{E_N}}$ is bounded uniformly in $N$. 
Moreover,  since $S_N u \to u$ as $N \to \infty$,  we deduce 
$\psi_{|_{E_N}} (u)  \to  \psi(u)$ for any $u\in L^p(\R)$. 
Thus, by Lebesgue's dominated convergence theorem, 
$\big|\psi_{|_{E_N}} -\psi\big|_{L^q(L^p, d\rho)}$ tends to $0$.
\hfill\qed
\medskip

\proof(of Proposition \ref{prop:prop6}) 
Let $\phi \in C_b(L^p(\R); \R)$. Then,  
$\phi$ is continuous and bounded on $E_N$ since $E_N \subset L^p(\R).$ 
Thus, if we choose $\phi_N:=\phi_{|_{E_N}}$, with the notation of Lemma \ref{lem:cylindrical_approx},
then $\tilde{P}_t^N \phi_N (y)= \E(\phi_N(X_N (t,y)))$ is well-defined for $y\in E_N$.
Hence
\begin{eqnarray*}
\int_{E_N} |\tilde{P}_t^N \phi_N (y)-\bar{\phi}_N^{\tilde{\rho}_N}|^2 d\tilde{\rho}_N(y) 
&=& 
\int_{E_N}   |\tilde{P}_t^N \phi_N (\Pi_N y)-\bar{\phi}_N^{\tilde{\rho}_N}|^2 
d\tilde{\rho}_N(\Pi_N y) \int_{E_N^{\perp}} d\mu_N^{\perp}((I-\Pi_N)y) \\
&=& \int_{L^p}  |\tilde{P}_t^N \phi_N (\Pi_N y')-\bar{\phi}_N^{\tilde{\rho}_N}|^2 
d\rho_N(y') \\
&=& \int_{L^p}  |\tilde{P}_t^N \phi_N (y')-\bar{\phi}_N^{\tilde{\rho}_N}|^2 
d{\rho}_N(y') 
\end{eqnarray*}
where we have used the following equalities: 
\begin{equation*}
\phi_N(\Pi_N y') =\phi_{|_{E_N}}(\Pi_N y')=
\phi (S_N \Pi_N y')=\phi(S_N y')=\phi_{|_{E_N}}(y')=\phi_N(y'). 
\end{equation*}
Thus, it follows that 
\begin{eqnarray*}
&&\Big|\int_{E_N} |\tilde{P}_t^N \phi_N(y)-\bar{\phi}_N^{\tilde{\rho}_N}|^2  d\tilde{\rho}_N (y) 
-\int_{L^p} |P_t \phi(y) -\bar{\phi}^{\rho}|^2  d\rho(y) \Big| \\
&=& 
\Big|\int_{L^p} |\tilde{P}_t^N \phi_N(y)-\bar{\phi}_N^{\tilde{\rho}_N}|^2  d{\rho}_N (y) 
-\int_{L^p} |P_t \phi(y) -\bar{\phi}^{\rho}|^2  d\rho(y) \Big| \\
&& \le 
\int_{L^p}|\tilde{P}_t^N \phi_N(y)-P_t \phi(y)|^2 d\rho_N(y) \\
&& +
\Big|\int_{L^p}|P_t\phi(y)-\bar{\phi}^{\rho}|^2 d\rho_N(y)-\int_{L^p}|P_t \phi(y)-\bar{\phi}^{\rho}|^2 d\rho(y)\Big|
+\Big|\bar{\phi}^{\rho}-\bar{\phi}_N^{\tilde{\rho}_N} \Big|^2 .
\end{eqnarray*}
Similar arguments as in the proof of  \eqref{invariancerho}
imply that all terms in the right hand side converge to zero as $N \to +\infty$. Indeed, note that 
$|P_t \phi(y)-\bar{\phi}^{\rho}|^2 \le 4|\phi|_{L^{\infty}}^2$ ;
thus, the second term on the right hand side of the last inequality  is estimated by 
\begin{equation*}
 4\Gamma |\Gamma_N^{-1}-\Gamma^{-1}| |\phi|_{L^{\infty}}^2 
+ 4|\phi|_{L^{\infty}}^2 \Gamma_N^{-1} \int_{L^p}\big| e^{-\frac{1}{4}|S_N y|_{L^4}^4}
-e^{-\frac{1}{4}|y|_{L^4}^4} \big| d\mu(y), 
\end{equation*}
and this quantity tends to zero by (\ref{bound:normalizedcst}) and the dominated convergence Theorem. 
The last term is estimated as follows, using the same computation as above.
\begin{eqnarray*}
\Big|\bar{\phi}^{\rho}-\bar{\phi}_N^{\tilde{\rho}_N} \Big| 
&&\le \Big|\int_{L^p}(\phi(y)-\phi_N(y)) d\rho(y) \Big| \\
&& + 
\Gamma |\Gamma_N^{-1}-\Gamma^{-1}|  |\phi|_{L^{\infty}} + |\phi|_{L^{\infty}} \Gamma_N^{-1}
\int_{L^p}\big|e^{-\frac{1}{4}|S_N y|_{L^4}^4}
-e^{-\frac{1}{4}|y|_{L^4}^4} \big|d\mu(y),
\end{eqnarray*}
which tends to zero when $N\to +\infty$, applying Lemma \ref{lem:cylindrical_approx} with $q=1$.   
Lastly, the choice of approximation
$X_N$ defined in Section 5 allows us to apply Proposition \ref{prop:prop4}  
to conclude that $X_N(\cdot, y) \to X(\cdot,y)$ in $C([0, T], L^p(\R))$ with $p\ge 3$ a.s., for any $y\in \bar{\mathcal{O}}$.
Therefore, the first term 
\begin{equation*}
\int_{L^p}|\tilde{P}_t^N \phi_N(y)-P_t\phi(y)|^2 d\rho_N(y) 
\le \int_{\bar{\mathcal{O}}} |\E(\phi(X_N(t, S_N y))-\phi(X(t,y)))|^2 d\mu (y)
\end{equation*}
converges to zero as $N\to +\infty$ by Lebesgue's dominated convergence. 
\hfill\qed
\medskip

\begin{rem}
Using the argument in \cite{odasso}, we could use the strong Feller property (Proposition \ref{p8}) and a coupling argument to prove that 
exponential convergence to equilibrium holds for any initial data. This may seem better than the result of Theorem \ref{thm:thm2} but,
contrary to Theorem \ref{thm:thm2}, the convergence rate given in the proof with such a coupling argument is difficult to write explicitly 
and  is very small.
\end{rem}

\begin{rem}
\label{rem7.2}
Note that if we consider, as in \cite{bs}, the same problem but posed on a bounded interval $D$ of $\R$, without quadratic potential, that is
$$
dX = (i-\gamma) (-\partial_x^2 X +\lambda |X|^2 X) dt +\sqrt{2\gamma} dW,
$$
then the corresponding Gibbs measure is supported in $L^2(D)$ (see \cite{bs}). All the arguments above may obviously be adapted to this
case if $\lambda >0$. Moreover, in the focusing case $\lambda=-1$, one may proceed as in \cite{cfl}, with $p=4$, $r=3$ and $\sigma=Id$, 
that is considering the modified dynamics
\begin{equation}
\label{foc}
dX= (i-\gamma) (-\partial_x^2 X - |X|^2 X +6 \kappa |X|_{L^2}^4 X) dt +\sqrt{2\gamma} dW
\end{equation}
corresponding to the Hamiltonian
$$
S(X) = \frac12 |\partial_x X|^2_{L^2} -\frac14 |X|_{L^4}^4 +\kappa |X|_{L^2}^6.
$$
For the associated purely dissipative dynamics
\begin{equation}
\label{diss-foc}
dX=\gamma (\partial_x^2 X + |X|^2 X - 6 \kappa |X|_{L^2}^4 X) dt + \sqrt{2\gamma} dW,
\end{equation}
it is not difficult to check that $\eta^h(t)=D_{X_0} X(t,X_0) h$ formally satisfies, for $\kappa$ large enough, an
estimate similar to \eqref{estimeta}, which allows to prove the Poincar\'e inequality for finite dimensional approximations
of the Gibbs measure, uniformly in the approximation parameter. All the other arguments seem to work in this case, and
this indicates that we easily recover with our method the convergence result of \cite{cfl}, at least in the case $p=4$.
Now, in our situation, where the domain is the whole space, with the addition of the quadratic potential $V(x)=x^2$, 
one has to take account of the fact that $\mu(L^2)=0$ (where we recall that $\mu$ is the Gaussian measure defined 
by the quadratic part of the Hamiltonian) so that the natural Gibbs measure given by
$$
e^{-\kappa |X|_{L^2}^6+\frac14 |X|_{L^4}^4} \mu(dX)
$$
is actually equal to the trivial measure $\delta_0$.
\end{rem}

\section{The case of general chemical potential $\eta$}

In this section, we explain how to treat the case where the chemical potential $\eta$ in \eqref{eq:CGL} is positive, and possibly
larger than $\lambda_0=1$. 
We use the notations of Section 2, and in particular, the measure $\rho$ is defined by \eqref{gibbs}.
It is clear that Proposition \ref{invariance} applies also to such a Gibbs measure through a finite dimensional approximation 
(see the arguments below) and we infer the invariance of the measure under the flow given by (\ref{eq:CGL}). 
Now, in order to prove the Poincar\'e inequality (\ref{ineq:poincare}) in this case, 
we write $S(u)$ in the following way :
\begin{equation*}
S(u)= \frac{1}{2}\int_{\R} |(-H)^{1/2}u|^2 dx + V_1(u)+V_2(u) 
\end{equation*}
with 
\begin{equation*}
V_1(u)=\int_{\R} F_1(x,|u(x)|^2) dx, \quad V_2(u)=\int_{\R} F_2(x, |u(x)|^2) dx,  
\end{equation*}
where $F_1$ and $F_2$ are defined by 
\begin{equation*}
F_1(x,y):=\left\{
\begin{array}{ll}
&\frac{1}{4}y^2 -\frac{\eta}{2}y,\quad y \ge \eta, \\[0.35cm]
&\Theta(x) \Big(\frac{1}{4} y^2 - \frac{\eta}{2} y\Big)+(1-\Theta(x))\Big(-\frac{\eta^2}{4}\Big), \quad 0\le y  <\eta,
\end{array}
\right.\\
\end{equation*}
and 
\begin{eqnarray*}
F_2(x,y)&:=& \frac14 y^2 -\frac{\eta}{2}y-F_1(x,y)\\
& =& (1-\Theta(x))\Big(\frac14 y^2-\frac{\eta}{2}y+\frac{\eta^2}{4}\Big)\un_{0\le y <\eta}.
\end{eqnarray*}
Here, $\Theta(x)$ is a cut-off function satisfying 
\begin{equation*}
\Theta \in C^{\infty}(\R), \quad 0 \le \Theta(x) \le 1, 
\quad \Theta(x)=1~ \mbox{for}~ |x| \ge \sqrt{2\eta}, \quad \Theta(x)=0~ \mbox{for}~ |x| \le \sqrt{3\eta/2}. 
\end{equation*}
Note that $F_1(x,y)$ is a convex function of $y\in \R^+$, for all $x\in \R$.
With these definitions, the Gibbs measure is rewritten under the following form :
\begin{equation*}
 \rho(du)= \Gamma^{-1} e^{-V_1(u)} e^{-V_2(u)} \mu(du), \quad \mathrm{with} \quad
\Gamma= \int_{L^p} e^{-V_1(u)} e^{-V_2(u)} d\mu(u),
\end{equation*}
and $\mu$ is the Gaussian measure defined in Section 2.
Note that $V_2(u)$ satisfies 
$e^{-C_\eta} \le e^{-V_2(u)} \le 1$ for some constant $C_\eta>0$ depending only on $\eta$, so that $e^{-V_2(u)}\mu(du)$ 
makes sense. Indeed, 
\begin{eqnarray*}
0\le V_2(u)&=& \int_{\R} F_2(x, |u|^2) dx= \int_{\R} (1-\Theta(x))\Big(\frac{1}{4}|u|^4 -\frac{\eta}{2}|u|^2
+\frac{\eta^2}{4}\Big) \un_{|u|^2 <\eta} dx \\
&\le& \frac{\eta^2}{2}\int_{\R} (1-\Theta(x)) dx \le \int_{|x| \le \sqrt{2\eta}} \frac{\eta^2}{2} dx=\sqrt{2} \eta^{\frac{5}{2}}. 
\end{eqnarray*}
Note that the finiteness of $\rho$ implies $e^{-V_1(u)}$ is integrable with respect to $ e^{-V_2(u)} d\mu$; then the boundedness
of $V_2$ implies that the measure $e^{-V_1(u)} \mu(du)$ is also finite.
\vspace{3mm}

Let us consider, for a large $N \in \N$ satisfying $\frac{2N_0+1}{2N+1} \le \frac{1}{4}$,  the measures defined on $E_N$ :
\begin{equation*}
\tilde{\rho}'_N(du):=(\tilde{\Gamma}_N')^{-1}e^{-V_2(S_N u)}e^{-V_1(S_N u)} \mu_N(du), 
\quad \mathrm{ and} \quad \pi_N(du):= e^{-V_1(S_N u)} \mu_N(du), 
\end{equation*}
where $\tilde{\Gamma}_N'= \int_{E_N} e^{-V_1(S_N u)} e^{-V_2(S_N u)} d\mu_N(u)$, and $\mu_N$ is defined in Section 2.
We will see that $\pi_N$ satisfies the Poincar\'e inequality, i.e. for any $\phi \in C^1_b(E_N, \R)$
\begin{equation} \label{ineq:poincare1}
\int_{E_N} |D\phi(y)|_{E_N}^2 d\pi_N(y) \ge \int_{E_N} |\phi(y)-\bar{\phi}^{\pi_N}|^2 d\pi_N(y),
\end{equation}
which implies a Poincar\'e inequality for $\tilde{\rho}'_N$. Indeed, assuming that \eqref{ineq:poincare1} holds,
\begin{equation*}
\begin{array}{l}
\ds e^{C_\eta}\int_{E_N} |D\phi(y)|_{E_N}^2 d\tilde{\rho}'_N(y)
= e^{C_\eta} (\tilde{\Gamma}_N')^{-1} \int_{E_N} |D\phi(y)|_{E_N}^2 e^{-V_2(S_N y)} d\pi_N(y)\\[0.4cm]
\ds \ge  (\tilde{\Gamma}_N')^{-1} \int_{E_N}|D\phi(y)|_{E_N}^2  d\pi_N(y)
\ge   (\tilde{\Gamma}_N')^{-1}  \int_{E_N} |\phi(y)-\bar{\phi}^{\pi_N}|^2 d\pi_N(y)\\[0.4cm]
\ds \ge(\tilde{\Gamma}_N')^{-1}\int_{E_N} |\phi(y)-\bar{\phi}^{\pi_N}|^2  e^{-V_2(S_N y)}  d\pi_N(y)
=  \int_{E_N} |\phi(y)-\bar{\phi}^{\pi_N}|^2 d\tilde{\rho}'_N(y),\\
\end{array}
\end{equation*}
and the Poincar\'e inequality for $\tilde{\rho}'_N$ follows from the obvious fact that 
$$
\int_{E_N} |\phi(y)-\bar{\phi}^{\tilde{\rho}'_N}|^2 d\tilde{\rho}'_N(y) =\inf_{C\in \R} \int_{E_N} |\phi(y)-C|^2 d\tilde{\rho}'_N(y).
$$
For (\ref{ineq:poincare1}), the following two points are essential according to Section 4. 
The first point is the dissipative inequality (\ref{ineq:dissipativity}), with Equation \eqref{equ:dissipative}
replaced by the equation 
\begin{equation*}
dY=\gamma(HY-2 S_N(\p_y F_1(x, |S_N Y|^2) S_N Y))dt +\sqrt{2\gamma} \Pi_N dW,
\end{equation*}
and the second point is the exponential decay of $\eta_N^h:= D_y Y(t,y).h$.
One may check for all $w,z \in E_N$ and $\gamma>0$,
\begin{eqnarray} \label{diss}
&&(\gamma H(w-z) -2\gamma S_N(\p_y F_1(\cdot,|S_N w|^2)S_N w -\p_y F_1(\cdot, |S_N z|^2)S_N z), w-z)_{E_N} \\
&=& -\gamma |(-H)^{1/2}(w-z)|_{L^2}^2 
-\gamma \mathrm{Re} \int_{\R} \left\{\p_y F_1(x, |S_N w|^2) + \p_y F_1(x, |S_N z|^2)\right\}|S_N w-S_N z|^2 dx \nonumber\\
&&-\gamma \mathrm{Re} \int_{\R} \left\{\p_y F_1(x, |S_N w|^2)-\p_y F_1(x, |S_N z|^2)\right\}(S_N w+S_N z)(\overline{S_N w} -\overline{S_N z}) dx. \nonumber
\end{eqnarray}
The second term on the right hand side above can be estimated as follows. 
\begin{eqnarray*}
&&-\gamma \mathrm{Re} \int_{\R} \p_y F_1(x, |S_N w|^2)|S_N w-S_N z|^2 dx \\ 
&=& -\gamma \mathrm{Re} \int_{\R}  \un_{|S_N w|^2 \ge \eta} \Big(\frac12 |S_N w|^2 -\frac{\eta}{2}\Big)|S_N w-S_N z|^2 dx \\
&& \hspace{3mm}-\gamma \mathrm{Re} \int_{\R}  \un_{|S_N w|^2 < \eta} \Theta(x) \Big(\frac12 |S_N w|^2 -\frac{\eta}{2}\Big)| S_N w-S_N z|^2 dx. \\
\end{eqnarray*}
Note that the first term on the right hand side above is non positive while the second term is bounded by
\begin{eqnarray*}
\frac{\gamma \eta}{2} \int_{\R} \Theta(x) |S_N w -S_N z|^2 dx 
& \le&\frac{\gamma}{2} \int_{|x|\ge \sqrt{\frac{3\eta}{2}}} \eta |S_N w-S_N z|^2 dx\\
&\le & \frac{\gamma}{3} \int_{\R} |x|^2 |(S_Nw-S_Nz)(x)|^2 dx\le
\frac{\gamma}{3} \int_{\R} |(-H)^{1/2}(S_Nw-S_N z)|^2 dx \\
&\le& \frac{\gamma}{3}  |(-H)^{1/2}(w -z)|_{L^2}^2.
\end{eqnarray*}
Thus 
\begin{equation*}
-\gamma \mathrm{Re} \int_{\R} \left\{\p_y F_1(x, |S_N w|^2) + \p_y F_1(x, |S_N z|^2)\right\}|S_N w-S_N z|^2 dx 
\le \frac{2\gamma}{3}  |(-H)^{1/2}(w -z)|_{L^2}^2,
\end{equation*}
and this term can be absorbed by the first term on the right hand side of \eqref{diss}. 
On the other hand, the third term, which is equal to
$$
-\gamma \mathrm{Re} \int_{\R} \left\{\p_y F_1(x, |S_N w|^2)-\p_y F_1(x, |S_N z|^2)\right\}(|S_N w|^2 -|S_N z|^2) dx
$$
is clearly non positive, thanks to the convexity of $F_1(x,\cdot)$.
Hence, the dissipative inequality holds. 

Consider now $\eta_N^h:= D_y Y(t,y).h$ which is the solution of
\begin{equation*}
\left\{
\begin{array}{ll}
\p_t \eta &= \gamma H \eta -2\gamma S_N \Big[\p_y F_1 (x, |S_N Y|^2)S_N \eta 
+ 2\p_y^2 F_1(x, |S_N Y|^2)\mathrm{Re}(\overline{S_N Y} S_N \eta) S_N Y \Big] \\ 
\eta(0) &= h \in E_N.
\end{array}
 \right.
\end{equation*}
Similarly as above, taking the product with $\eta_N^h$,  
\begin{eqnarray*}
\frac{1}{2} \frac{d}{dt}|\eta_N^h|_{E_N}^2 + \gamma(-H \eta_N^h, \eta_N^h)_{E_N} 
&=& 
-2\gamma \mathrm{Re} \int_{\R} \p_y F_1(x, |S_N Y|^2)|S_N \eta_N^h|^2 dx \\
&& -4\gamma \int_{\R} \p_y^2 F_1(x, |S_N Y|^2) \left\{\mathrm{Re}(\overline{S_N Y} S_N \eta_N^h) \right\}^2 dx \\
&& \le \frac{2\gamma}{3}  \int_{\R} |(-H)^{1/2} \eta_N^h|^2 dx,
\end{eqnarray*}
and the exponential decay of $\eta_N^h$ follows. Thus, all the results of Section 4 also hold in the case $\eta>0$.
\vspace{3mm}

Next, the results of Sections 5 and 6 are proved exactly in the same way as in the case $\eta=0$. It remains to obtain 
the convergence of the measure $\tilde{\rho}'_N$ (Section 7); for that purpose, 
consider $N$ large enough such that $\frac{2N_0+1}{2N+1} \le \frac{1}{4}$, define $A_N:=-H-\eta(S_N-\Pi_{N_0})$ and 
consider the associated Gaussian measure $\tilde{\mu}_N$ on $E_N$ which, up to a normalizing constant, is equal to
$e^{-\frac{1}{2}(A_N u, u)_{E_N}} du$, and converges to
$\tilde{\mu}_{\eta}$ as $N \to +\infty$. 
Note that we may write 
\begin{equation*}
\tilde{\rho}'_N(du)=(\tilde{\Gamma}_N)^{-1} e^{-\tilde{V}(S_N u)} d\tilde{\mu}_N(du),  
\end{equation*}
with $\tilde{\Gamma}_N= \int_{E_N} e^{-\tilde{V}(S_N u)} d\tilde{\mu}_N(du)$. 
We have $\tilde{\rho}'_N \otimes \tilde{\mu}_N^{\perp} \to \rho$ as $N \to +\infty$, for $p>2$, provided 
\begin{equation*}
\lim_{NÊ\to +\infty} \tilde{\mu}_{\eta}(u\in L^p(\R),~ |\tilde{V}(S_N u)-\tilde{V}(u)| >\delta)=0, 
\end{equation*}
for any $\delta>0$.
This may be shown similarly as in Lemma 3.3 in \cite{btt}.
\vspace{3mm}

\section{Appendix}

Here we give a proof of Lemma 2.1.
\vspace{3mm}

\noindent
{\bf Proof of Lemma 2.1.} 
We apply the Kolmogorov test (see \cite{dpz}) in order to investigate the regularity of $Z_{\infty}(t)$. 
For $t,s \in [-T',T],$ with $s<t$,
\begin{eqnarray*}
Z_{\infty}(t,x)-Z_{\infty}(s,x)&=& \sqrt{2\gamma} \sum_{k\in \N} \Big[\int_s^t e^{-\lambda_k^2 (i+\gamma)(t-\sigma)} (d\beta_k^R (\sigma) h_k(x)+id \beta_k^I(\sigma) h_k(x)) \\
             &+& \int_{-\infty}^s (e^{-\lambda_k^2(i+\gamma)(t-\sigma)}-e^{-\lambda_k^2(i+\gamma)(s-\sigma)}) 
             (d\beta_k^R (\sigma) h_k(x)+id \beta_k^I(\sigma) h_k(x)) \Big].
\end{eqnarray*}
We set 
\begin{eqnarray*}
f_{1,t,s}(\omega,x)&:=& \sum_{k} \int_s^t e^{-\lambda_k^2 (i+\gamma)(t-\sigma)} (d\beta_k^R (\sigma) h_k(x)+id \beta_k^I(\sigma) h_k(x)),\\ 
f_{2,t,s}(\omega,x)&:=& \sum_{k} \int_{-\infty}^s (e^{-\lambda_k^2(i+\gamma)(t-\sigma)}-e^{-\lambda_k^2(i+\gamma)(s-\sigma)}) 
(d\beta_k^R (\sigma) h_k(x)+id \beta_k^I(\sigma) h_k(x)),
\end{eqnarray*} 
and we will make use of the Minkowski inequality; for $q \ge p,$ 
$$|f_{j,t,s}(\omega,x)|_{L^q_{\omega} L^p_{x}} \le |f_{j,t,s}(\omega,x)|_{L^p_x L^q_{\omega}}, \quad j=1,2, \quad t,s\in [0,T].$$
We calculate first  
\begin{eqnarray*}
& & \E\Big(\Big|\sum_{k} \int_s^t e^{-\lambda_k^2 (i+\gamma)(t-\sigma)}(d\beta_k^R (\sigma) h_k(x)+id \beta_k^I(\sigma) h_k(x)) \Big|^2 \Big) \\
&\le& 2 \E \sum_{k} \int_s^t e^{-2 \gamma \lambda_k^2 (t-\sigma)} d\sigma |h_k(x)|^2 \\ 
&\le& \sum_{k} \frac{1}{\gamma \lambda_k^2} (1-e^{-2\gamma \lambda_k^2 (t-s)})|h_k(x)|^2 \\
&\le & C_{\alpha, \gamma} \sum_k \lambda_k^{2(\alpha-1)}|t-s|^{\alpha} |h_k(x)|^2,
\end{eqnarray*}
for any $\alpha \in [0,1]$, and $\gamma >0$. Since moreover $\{f_{1,t,s}(\omega, x)\}_{t \in [0,T]}$ is a Gaussian process, we deduce that  
for any $m \in \N \setminus\{0\}$, 
\begin{eqnarray*}
&&\E\Big(\Big|\sum_{k} \int_s^t e^{-\lambda_k^2 (i+\gamma)(t-\sigma)}(d\beta_k^R (\sigma) h_k(x)+id \beta_k^I(\sigma) h_k(x)) \Big|^{2m}\Big)\\
&&\le (C_{\alpha, \gamma})^m |t-s|^{\alpha m} \Big(\sum_k \lambda_k^{2(\alpha-1)} |h_k(x)|^2 \Big)^m.
\end{eqnarray*}
Therefore, for $q=2m \ge p$, and $\alpha\in [0,1],$ 
\begin{eqnarray*}
\E \Big(|f_{1,t,s}(x)|_{L^p(\R)}^{2m}\Big)^{\frac{1}{2m}} 
&=&|f_{1,t,s}|_{L^{2m}_{\omega} L^p_x} \\  
&\le& |f_{1,t,s}|_{L^p_x L^{2m}_{\omega}} \\
&\le& (C_{\alpha, \gamma})^{1/2} |t-s|^{\alpha/2} \Big|\Big(\sum_{k} \lambda_k^{2(\alpha-1)} (h_k(x))^2\Big)^{\frac{1}{2}}\Big|_{L^p_x}.
\end{eqnarray*}
On the other hand, 
\begin{eqnarray*}
\Big|\Big(\sum_k \lambda_k^{2(\alpha-1)} |h_k(x)|^2 \Big)^{1/2}\Big|_{L^p_x} 
= \Big|\sum_k \lambda_k^{2(\alpha-1)} |h_k(x)|^2 \Big|_{L^{p/2}_x}^{1/2}  \le  \Big(\sum_k \lambda_k^{2(\alpha-1)} |h_k|_{L^p_x}^2\Big)^{1/2}.
\end{eqnarray*}
Here we remark that for all $p\ge 4,$ there exists $C_p>0$ such that
$|h_k|_{L^p(\R)} \le C_p \lambda_k^{-\frac{1}{6}}$ 
(see \cite{yz}, Lemma 3.2), and by interpolation, if $2 \le p\le 4,$ 
$|h_k|_{L^p(\R)} \le C_p \lambda_k^{-\frac{1}{6}(2-\frac{4}{p})}$. 
Therefore, using the notation $\theta(p)$ in (\ref{def:thetap}),   
\begin{equation*}
|\Big(\sum_k \lambda_k^{2(\alpha-1)} |h_k(x)|^2 \Big)^{1/2}|_{L^p_x} 
\le C_p \Big(\sum_k \lambda_k^{2(\alpha-1-\frac{1}{6}\theta(p))} \Big)^{1/2} 
\end{equation*}
and the series converges if $p>2$ and $\alpha<\theta(p)/6.$ Namely, 
\begin{equation} \label{est:f_1}
\E(|f_{1,t,s}|_{L^p_x}^{2m}) \le C_{p, \alpha, \gamma, m} |t-s|^{\alpha m} 
\end{equation}
for $\alpha<\theta(p)/6$ and $p>2$. 
Similarly, for the second term $f_{2,t,s}(\omega,x)$, 
\begin{eqnarray*}
&& \E(|\sum_k \int_{-\infty}^s (e^{-\lambda_k^2(i+\gamma)(t-\sigma)}-e^{-\lambda_k^2(i+\gamma)(s-\sigma)}) (d\beta_k^R (\sigma) h_k(x)+id \beta_k^I(\sigma) h_k(x))|^2) \\
&\le& 2\E \sum_k  \int_{-\infty}^s |e^{-\lambda_k^2(i+\gamma)(t-\sigma)}-e^{-\lambda_k^2(i+\gamma)(s-\sigma)}|^2 d\sigma |h_k(x)|^2 \\ 
&\le & 2\sum_k \int_{-\infty}^s e^{-2\lambda_k^2 \gamma (s-\sigma)} |e^{-\lambda_k^2 (i+\gamma)(t-s)}-1|^2 d\sigma |h_k(x)|^2\\
&\le & \sum_k \frac{1}{\gamma \lambda_k^2} |e^{-\lambda_k^2(i+\gamma)(t-s)}-1|^2 |h_k(x)|^2.
\end{eqnarray*}
Here, we note that 
$$
|e^{-\lambda_k^2(i+\gamma)(t-s)}-1| 
\le C_{\alpha, \gamma} \lambda_k^{2\alpha} |t-s|^{\alpha},
$$
if $t-s \ge 0$ and $0 \le \alpha \le 1$.
Thus, 
\begin{eqnarray*}
&&\E(|\sum_k \int_{-\infty}^s (e^{-\lambda_k^2(i+\gamma)(t-\sigma)}-e^{-\lambda_k^2(i+\gamma)(s-\sigma)})(d\beta_k^R (\sigma) h_k(x)+id \beta_k^I(\sigma) h_k(x)) |^2)\\ 
&\le& C'_{\gamma, \alpha} \sum_k \lambda_k^{2(2\alpha-1)}|t-s|^{2\alpha} |h_k(x)|^2
\end{eqnarray*}
for all $\alpha \in [0,1],$ and we conclude with the same arguments 
as for the first term $f_{1,t,s}(\omega,x)$ that for $m\ge \frac{p}{2}>1$ and $\alpha <\theta(p)/12$, 
\begin{equation*}
\E(|f_{2,t,s}|_{L^p_x}^{2m}) \le C_{p,\alpha, \gamma, m} |t-s|^{2\alpha m}.
\end{equation*} 
We conclude by the Kolmogorov test
that $Z_{\infty}$ has a modification in $C^{\alpha'}([0,T], L^p(\R))$ for any $\alpha' <\theta(p)/12$ 
and $p>2$. In particular, 
\begin{equation*}
\E \Big(\sup_{t\in [0,T]}|Z_{\infty}(t)|_{L^p(\R)}\Big) \le M_{T,p}. 
\end{equation*}
Note that the same calculation with $L^p(\R)$ replaced by $\mathcal{W}^{s,p}(\R)$ can be performed, and shows that
$Z_{\infty}(t)$ has a modification in $C^{\alpha'}([0,T], \mathcal{W}^{s,p}(\R))$ for $p>2$, $0 \le s <\theta(p)/6$, and
$\alpha' <\theta(p)/12-s/2$.
Furthermore, the arguments above imply the same estimates for $\Pi_N Z_{\infty}$ and $(I-\Pi_N)Z_{\infty}$.
\hfill\qed 
\bigskip


\noindent 
{\bf Acknowledgement}
This work is supported by the JSPS KAKENHI Grant Numbers 15K04944 and 16KT0127, the ANR project 
ANR-12-MONU-0007 BECASIM, and the French government ``Investissements d'Avenir" program ANR-11-LABX-0020-01.
The authors are grateful to M. Kobayashi for helpful discussions. 

\end{document}